%% file: ribbonknotdb.tex
\documentclass[twoside]{article}
\input{header}

\theorems
\begin{document}
\papertitle{A Tabulation of Ribbon Knots in Tangle Form}
\paperauth{Andrey Boris Khesin}{University of Toronto}
\begin{paperabs}
It is known that there are 21 ribbon knots with 10 crossings or fewer.
We show that for every ribbon knot, there exists a tangle that satisfies two
properties associated with the knot.
First, under a specific closure, the closed tangle is equivalent to its
corresponding knot.
Second, under a different closure, the closed tangle is equivalent to the
unlink.
For each of these 21 ribbon knots, we present a 4-strand tangle that
satisfies these properties.
We provide diagrams of these tangles and also express them in planar diagram
notation.
\end{paperabs}
\begin{paper}
\papersec{Introduction}

There are 250 knots with 10 crossings or fewer.
The list of these knots is commonly referred to as the Rolfsen Table.
There are 21 knots in the Rolfsen Table, deemed \textit{ribbon knots}, that
satisfy a particular set of conditions.
These knots are $6_1$, $8_8$, $8_9$, $8_{20}$, $9_{27}$, $9_{41}$, $9_{46}$,
$10_3$, $10_{22}$, $10_{35}$, $10_{42}$, $10_{48}$, $10_{75}$, $10_{87}$,
$10_{99}$, $10_{123}$, $10_{129}$, $10_{137}$, $10_{140}$, $10_{153}$, and
$10_{155}$.
We will show that for each ribbon knot there exists a tangle that satisfies two
particular properties associated with the knot.
Furthermore, if such a tangle exists for a given knot, then that knot is a
ribbon knot.
Lastly, we will provide diagrams of such tangles for these 21 ribbon knots and
express those tangles in planar diagram notation.\\

\paperfig{Singularities}{\def\svgwidth{0.4\columnwidth}
{\centering\input{clasp.pdf_tex}}\hfill
\def\svgwidth{0.4\columnwidth}
{\centering\input{ribbon.pdf_tex}}\\

\noindent Clasp Singularity\hfill Ribbon Singularity}
{A clasp singularity and a ribbon singularity.
Let $M$ be the image of a disc $D$ that has been immersed into three-dimensional
space.
Every knot $K$ is the boundary of such an immersed disc $M$.
In three-dimensional space, $M$ can be drawn so that all of its
self-intersections are 1-manifolds.
Each of these self-intersections is a singularity curve.
Let $C$ be one of these curves.
We know that since the $C$ is a self-intersection of $M$, there are curves $C_1$
and $C_2$ in $D$ whose images are $C$.
Consider the endpoints of $C_1$ and $C_2$.
Two of these four endpoints have to be on the boundary of $D$ since the
endpoints of $C$ are on the boundary of $M$.
If either $C_1$ or $C_2$ has both of its endpoints on the boundary of $D$, the
singularity corresponding to the self-intersection $C$ is a ribbon singularity,
otherwise it is a clasp singularity.
The translucent regions are the immersed discs $M$, with the gray lines
representing the discs' boundaries, the knots $K$.
The black lines are the singularities, the curves $C$.}\\

We start with the definition of a ribbon knot.
Every knot is the boundary of a disc immersed into three-dimensional space.
Unless the immersion is trivial and its boundary is the unknot, this disc will
have self-intersections, or \textit{singularities}.

\paperdef{Ribbon}
{A \textit{ribbon knot} is a knot that can be drawn as the boundary of a disc
immersed into three-dimensional space that contains only ribbon singularities
(see \figSingularities).}

\papersvg{Immersion}{immersion}
{The first step in constructing a ribbon knot.
We start with a series of discs connected in a chain with arcs.
In this case, there are two such discs.
We then immerse our structure into three dimensional space, keeping the
immersion of the discs trivial.
The immersion of the arcs is non-trivial.
To turn this into a ribbon knot, all we have to do is unzip the image of the
arc into a pair of strands, a ribbon.
Unzipping the immersion shown creates the ribbon knot $6_1$.}\\

There is a way to construct a knot that is guaranteed to be a ribbon knot.
To do this, instead of immersing a single disc into three-dimensional space, we
immerse several such discs, each connected to the next by an arc, resulting in
what looks like a string of beads.
However, the immersion of each disc must be trivial.
We cannot let the discs have self-intersections as this might result in clasp
singularities.
Fortunately, there are no restrictions on the arcs, so they can freely pass
through the images of discs.
We present the image of the structure used to create the ribbon knot $6_1$ (see
\figImmersion).
Finally, after we immerse our structure we unzip each arc by turning it into a
\textit{ribbon}, two parallel strands tracing out the path of the arc.
Note that after the last step, the structure in the diagram above becomes the
knot in the diagram below.

\paperprp{Generation}
{The above process of generating a knot always produces a ribbon knot.}
\begin{proof}[Proof of Proposition]
After applying the above process, the only singularities in the resulting knot
will arise where the ribbons intersect the embedded discs.
We see that the resulting curves of self-intersection do not touch the edges of
the discs.
However, since the curves go right down the width of the ribbons, they do touch
the boundaries of those ribbons.
Thus, all of the singularities in the knot will be ribbon singularities since
the curves of self-intersection have both endpoints on the ribbon, as opposed to
having one on the ribbon and the other on the disc's boundary.
Thus the knot will be a ribbon knot by \defRibbon (see \figSingularities).
\end{proof}

\papersvg{Presentation}{presentation}{A ribbon presentation of $6_1$.
Note that there are two rectangular components which are the images of the
boundaries of two discs.
The ribbon connecting them is constructed by unzipping the image of an arc that
connected the discs, thereby replacing a single strand with a pair.
Since the only singularities in ribbon presentations can arise where the ribbons
go through the components, any ribbon presentation is a knot diagram of a ribbon
knot.}\\

A \textit{ribbon presentation} is a knot diagram which consists of components
and ribbons, all satisfying certain properties.
If we represent the components as a row of several large rectangles, and the
ribbons as pairs of parallel strands which connect each component to the next,
then the ribbons may go through and loop around the components, but are not
allowed to go through other ribbons.
If a knot diagram satisfies these properties, it is a called a ribbon
presentation.
We give a ribbon presentation of knot $6_1$ in \figPresentation.

Note that a ribbon knot generated by the process described earlier can be drawn
as a ribbon presentation.
It suffices to map each disc to a rectangle and when we unzip the arcs to form
ribbons, the result will become a ribbon presentation.

Any ribbon presentation is a knot diagram of a ribbon knot.
If we have a ribbon presentation of some knot, we start by filling in the image
of the disc whose boundary is the knot.
We see that all of the singularities occur where the ribbons pass through the
large rectangular components.
We also observe that all of those singularities are ribbon singularities so any
knot expressed as a ribbon presentation is a ribbon knot, making the name
appropriate (see \figSingularities).
Note that while this is similar to the proof of \prpGeneration, we emphasize
that it is not necessary to use the process of generating ribbon knots described
earlier to find out if a ribbon presentation is a knot diagram of a ribbon knot.

\papersec{Tangles}

Another important structure that we need to define is that of a tangle.
A \textit{tangle} is a collection of crossings with a certain number of loose
ends.
An $n$-strand tangle will have $2n$ such ends.\\

\paperfig{Tangle}{\size{9}{
\hspace{0.19\columnwidth}$t_1$\hspace{0.17\columnwidth}$t_2$
\hspace{0.05\columnwidth}\dots\hspace{0.04\columnwidth}$t_{n-1}$
\hspace{0.14\columnwidth}$t_n$\\
\svgc{tangle}\\
\indent\hspace{0.19\columnwidth}$b_1$\hspace{0.16\columnwidth}$b_2$
\hspace{0.04\columnwidth}\dots\hspace{0.04\columnwidth}$b_{n-1}$
\hspace{0.13\columnwidth}$b_n$}}
{An $n$-strand tangle T.
In a tangle, the ends do not have to be drawn along the top and bottom of the
tangle, this is done for convenience.
The $2n$ ends can be located anywhere on the perimeter.
We label the ends of an $n$-strand tangle from $t_1$ to $t_n$ along the top and
$b_1$ to $b_n$ along the bottom.
Here, the dashed strands are tangled and cross many times, but link up the
corresponding ends along the top and the bottom the way they do in a pure
braid.}\\

For convenience, a tangle is often drawn like a pure braid, with $n$ ends along
both the top and the bottom.
Thus, if we were to follow any strand of a tangle, we would end right below
from where we started.
We denote the tops of the strands by $t_1,~t_2,~\dots,~t_n$, and the
bottoms by $b_1,~b_2,~\dots,~b_n$ (see \figTangle).
This means that $t_k$ is connected through the tangle to $b_k$, where
$k=1,~\dots,~n$ is the strand index.
We will also denote the space of all $n$-strand tangles with $\mathcal{T}_n$.
Note that the space of all knots is isomorphic to $\mathcal{T}_1$ as a knot is
merely a 1-tangle with $t_1$ joined to $b_1$.

A \textit{closure} of a tangle is an operation which stitches some of the ends
of the tangle together.
Note that ends can only be stitched together in pairs.
Unless we are applying a closure to a tangle in $\mathcal{T}_1$, stitching
together a pair of ends of a tangle in $\mathcal{T}_n$ results in a tangle in
$\mathcal{T}_{n-1}$.
By appropriately stitching every single end of a tangle, we can turn that tangle
into a knot.\\

\paperfig{Top}{\size{9}{
\hspace{0.1\columnwidth}$t_1$\hspace{0.07\columnwidth}$t_2$
\hspace{0.07\columnwidth}$t_3$\hspace{0.14\columnwidth}\dots
\hspace{0.08\columnwidth}$t_{2n-2}$\hspace{0.04\columnwidth}$t_{2n-1}$
\hspace{0.04\columnwidth}$t_{2n}$\\
\svgc{top}\\
\indent\hspace{0.1\columnwidth}$b_1$\hspace{0.07\columnwidth}$b_2$
\hspace{0.06\columnwidth}$b_3$\hspace{0.125\columnwidth}\dots
\hspace{0.08\columnwidth}$b_{2n-2}$\hspace{0.03\columnwidth}$b_{2n-1}$
\hspace{0.03\columnwidth}$b_{2n}$}}
{The top closure of a $2n$-strand tangle T.
The top closure leaves every bottom end open and leaves all the upper ends
closed.
After the closure, the tangle has only $n$ strands instead of $2n$.
We denote this closure $\tau(T)$.}

\paperdef{Top}{For a $2n$-strand tangle $T$, we define the
\textit{top closure} $\tau:\mathcal{T}_{2n}\to\mathcal{T}_n$ of $T$ to be the
closure that stitches together the pairs
$(t_1,~t_2),~(t_3,~t_4),~\dots,~(t_{2n-1},~t_{2n})$ to create the $n$-strand
tangle $\tau(T)$ (see \figTop).}

\paperdef{Bottom}{For a $2n$-strand tangle $T$, we define the
\textit{bottom closure} $\beta:\mathcal{T}_{2n}\to\mathcal{T}_n$ of $T$ to be
the closure that stitches together the pairs $(b_1,~b_2),~(b_3,~b_4)$, \dots,
$(b_{2n-1},~b_{2n})$ to create the $n$-strand tangle $\beta(T)$.}

\paperfig{Full}{\size{9}{
\hspace{0.1\columnwidth}$t_1$\hspace{0.07\columnwidth}$t_2$
\hspace{0.07\columnwidth}$t_3$\hspace{0.14\columnwidth}\dots
\hspace{0.08\columnwidth}$t_{2n-2}$\hspace{0.04\columnwidth}$t_{2n-1}$
\hspace{0.04\columnwidth}$t_{2n}$\\
\svgc{full}\\
\indent\hspace{0.1\columnwidth}$b_1$\hspace{0.07\columnwidth}$b_2$
\hspace{0.06\columnwidth}$b_3$\hspace{0.125\columnwidth}\dots
\hspace{0.08\columnwidth}$b_{2n-2}$\hspace{0.03\columnwidth}$b_{2n-1}$
\hspace{0.03\columnwidth}$b_{2n}$}}
{The full closure of a $2n$-strand tangle T.
It is important to note that the full closure reproduces the top closure along
the bottom ends, and closes the top ends similarly, but offset by 1 strand.
This closure leaves ends $t_1$ and $t_{2n}$ open.
After the closure, the tangle has only 1 strand instead of $2n$.
We denote this closure $\phi(T)$.}

\paperdef{Full}{For a $2n$-strand tangle $T$, we define the
\textit{full closure} $\phi:\mathcal{T}_{2n}\to\mathcal{T}_1$ of $T$ to be the
closure that stitches together the pairs $(b_1,~b_2),~(t_2,~t_3),~(b_3,~b_4)$,
\dots, $(t_{2n-2},~t_{2n-1}),~(b_{2n-1},~b_{2n})$ to create the 1-strand tangle
$\phi(T)$ (see \figFull).}

Note that the stitchings applied by $\beta$ are a subset of those applied by
$\phi$.

If the inside of a $2n$-strand tangle $T$ is trivial, meaning that $T$ is an
\textit{untangle} and that it has no internal crossings, then closing the
strands along the top and bottom in adjacent pairs by applying both $\tau$ and
$\beta$ will allow us to create the $n$-component unlink, $U$.

\papereq{Unlink}{(\beta\circ\tau)(T)=U}{}
\begin{paperwhere}
\papervar{\beta}{the bottom closure}{}
\papervar{\tau}{the top closure}{}
\papervar{T}{a $2n$ strand untangle}{}
\papervar{U}{the $n$-component unlink}{}
\end{paperwhere}

We note two things about \eqUnlink.
First, the choice of order when applying $\beta$ and $\tau$ is irrelevant as the
operations commute.
Second, the relation \eqUnlink can still be satisfied even if $T$ is not an
untangle.
A simple example of such a tangle $T$ is a 2-strand tangle with just one
crossing.

Having defined our closures, we now prove our main result on the existence of a
common property of all ribbon knots.

\paperthm{Ribbon}
{A knot $K$ is a ribbon knot if and only if there exists a $2n$-strand tangle
$T$ such that $(\beta\circ\tau)(T)$ is the $n$-component unlink and $\phi(T)$
creates a 1-strand tangle equivalent to $K$.}
\begin{proof}
We begin by proving the ``only if'' direction.
First we show that there exist surfaces bounded by $(\beta\circ\tau)(T)$ that do
not have any singularities.

We assume that $(\beta\circ\tau)(T)$ is the $n$-component unlink.
Any $n$-component link is the immersion of $n$ discs into three-dimensional
space.
By applying a continuous spatial deformation $f$ to our unlink, we turn it back
into the boundaries of its original discs.
Applying the inverse transformation, $f^{\text-1}$, to the original discs, will
result in the boundaries of the discs being deformed to end up as our unlink,
while the surfaces will get continuously deformed.
Since the surfaces of the original discs did not intersect, neither will they
after deforming them continuously.
Thus, we can draw surfaces bounded by $(\beta\circ\tau)(T)$ with no
singularities.\\

\papersvg{Proof}{proof}
{A modified version of $(\beta\circ\tau)(T)$.
Here, we have connected each of the strands forming the stitchings along the top
to the next one with a pair of strands.
Note that this results in one long strand running along the top of the knot, so
these connections are just $\phi(T)$ with the added connection
$(t_1,~t_{2n})$.}\\

Since $(\beta\circ\tau)(T)$ is the $n$-component unlink, it is composed of $n$
unknots.
Each of those unknots can be connected to the next one with a thin ribbon by
adding a ribbon between neighbouring strands connecting the pairs
$(t_1,~t_2),~(t_3,~t_4),~\dots,~(t_{2n-1},t_{2n})$ (see \figProof).
As a result, we will have created $\phi(T)$ with the added stitching
$(t_1,~t_{2n})$.
Note that since all we added were thin ribbons, they only pass through the
surface whose boundary is the knot, thereby only forming ribbon singularities.
This is true for exactly the same reasons that a ribbon presentation is a knot
diagram of a ribbon knot.

As a result, we have shown that $(\beta\circ\tau)(T)$ contains no singularities
and the addition of our ribbons only added ribbon singularities.
Since the result is $\phi(T)$ with the stitching $(t_1,~t_{2n})$ and is
identical to our knot $K$, we have shown that $K$ is a ribbon knot, thereby
proving one direction of the theorem.

Conversely, we will show how, from any ribbon knot $K$, we can algorithmically
construct a tangle $T$ satisfying the necessary conditions.
First, consider an arbitrary $n$-strand tangle $T$.
We place all $2n$ of its ends along the top and number them from $t_1$ to
$t_{2n}$.
We will show that this can be used to construct a tangle with ends both along
the top and the bottom, which will later give us the desired tangle.\\

\paperfig{Lemma}{\def\svgwidth{0.4\columnwidth}
{\centering\input{lemmaone.pdf_tex}}\\

\vspace{-4.5em}\hspace{0.3in}\hspace{14.8ex}\size{20}{$\equiv$}\\

\vspace{-4.2em}\hfill\def\svgwidth{0.4\columnwidth}
{\centering\input{lemmatwo.pdf_tex}}\\}
{A visual representation of two tangles.
The equality shows that if there exists a tangle $T$, then there is also a
tangle $T'$ such that $\beta(T')=T$.}

\paperlem{Tangles}
{For any $n$-strand tangle $T$ with all $2n$ of its ends placed along the top,
there exists a $2n$-strand tangle $T'$ such that $\beta(T')=T$ (see \figLemma).}
\begin{proof}[Proof of Lemma]
Each of $T$'s $n$ strands starts at one of the ends along the top, comes down,
weaves through the tangle, and then comes out of the tangle through another end
along the top.
Thus, these $n$ strands effectively look like $n$ hanging loops that have been
twisted together.
All we need to do is to stretch each of these loops downwards until they are
lowered past the bottom edge of the tangle.
Then, we simply cut each of the loops, resulting in $2n$ strands, which we can
reorder so that their order matches that of their respective strands along the
upper edge.

What is left inside the boundaries of the tangle is a $2n$-strand tangle such
that applying the bottom closure to it gives $T$.
Thus, this tangle satisfies the properties of $T'$, which means that not only
does such a $T'$ always exist, but we also know how to construct it.
\end{proof}

\papersvg{Lowered}{lowered}
{A ribbon presentation of a knot after all of its intersections were enclosed in
one large tangle labeled $R$.
The components of the knots are depicted exactly as they are, while the
intersections of the ribbons and the components are not shown.}\\

We start our construction of $T$ from a ribbon presentation of the ribbon knot
$K$.
We enclose the lower halves of the components of that ribbon presentation of $K$
in one large tangle $R$.
Next, we move all of the ribbons so that the locations where they connect to the
components are above $R$, but all their crossings occur inside $R$
(see \figLowered).
This can be done since we can pull the ribbons into $R$, regardless of whether
they are twisted around the components or each other.\\

\papersvg{Twisted}{twisted}
{A ribbon presentation of a knot after all of its intersections were enclosed in
one large tangle labeled $R$ and all of the ribbons were pulled out of it.
The ribbons outside of the knots are depicted exactly as they are, while the
intersections of the components are not shown.}\\

After we have enclosed the knot in a tangle, we attempt to pull all of the
ribbons out of $R$.
To pull a given ribbon up out of $R$, we need to make sure that there is nothing
inside $R$ that is above the ribbon.
However, there will generally be parts of the components that the ribbon wraps
around.
What we do now is simply move those strands of the components up and around the
rest of the knot.
The components will end up back in $R$ and we will be free to raise the ribbons
up out of $R$ (see \figTwisted).
Since we have moved various other parts of the knot in and out of the tangle
$R$, we give the resulting tangle the name $R'$.
Note that at this point, the components in $R'$ are twisted, whereas the
components in $R$ only intersected the ribbons.

As this last step is the most complicated part of the proof, we show how this is
done for the simplest ribbon knot, $6_1$.
We start by examining the knot's ribbon presentation.\\

\papersvg{ExplainZero}{presentation}{The ribbon presentation of $6_1$.}\\

Our goal is to be able to lift the ribbon out from the knot so we position it
going nicely from left to right.
We start by moving the top of the ribbon to the left and the bottom of the
ribbon to the right.\\

\papersvg{ExplainOne}{explanationone}{A slightly deformed version of the ribbon
presentation of $6_1$.}

We now want to make the ribbon horizontal so we move the left half of the knot
down and the right half of the knot up.
This straightens the ribbon.\\

\papersvg{ExplainTwo}{explanationtwo}{The knot $6_1$ but deformed to make the
ribbon horizontal.}

There are two loops around the ribbon that are preventing us from pulling it up.
We remove the left one first by moving the back part of the loop backwards and
downwards and moving the front part towards us and downwards.
Although the component on the right is drawn smaller, this does not affect our
steps in deforming the knot.\\

\papersvg{ExplainThree}{explanationthree}{The knot $6_1$ with only one loop
preventing the ribbon from being pulled up.}

Finally, to remove the last loop, we slide it around the component on the right.
This leaves nothing above the ribbon, allowing us to pull it up above the rest
of the knot.
We then enclose everything below the ribbon in one large tangle.
This tangle is what we need.

\papersvg{ExplainFour}{explanationfour}{The knot $6_1$ with the ribbon pulled
up.}\\

To prove that the tangle in the diagram is the one we seek, we simply need to
apply our closures (see \figExplainFour).
If we were to apply the top and bottom closures to that tangle, we would get the
two-component unlink, as required.

Now, we continue the general proof, assuming that we have just pulled up the
ribbons of the knot we were working with.\\

\papersvg{Final}{final}
{Here, $T$ is a tangle such that $\beta(T)=R'$.
The existence of such a tangle is guaranteed by \lemTangles.
Note that this knot is $\phi(T)$ with the added stitching $(t_1,~t_{2n})$.}\\

Note that $R'$ is a tangle with an even number of ends placed along the top.
We will say that $R'$ has $2n$ ends.
By \lemTangles we know that there exists a tangle, which we will preliminarily
name $T$, such that $\beta(T)=R'$.
Note that the knot we have created with the tangle $R'$ at the bottom is
equivalent to $\phi(T)$ with the added stitching $(t_1,~t_{2n})$ since the top
part of the knot has no intersections anymore (see \figFinal).
Also note that this knot is still equivalent to $K$ as we have merely deformed
a ribbon presentation of $K$.
Thus we have shown that there exists a tangle $T$ satisfying one of the
conditions of the theorem.
All we have left to show is that $(\beta\circ\tau)(T)$ is the $n$-component
unlink.

We know that since every ribbon knot has a ribbon presentation, removing the
ribbons connecting the components of the knot results in the $n$-component
unlink (see \figPresentation).
This is true because when constructing a ribbon knot, we only immerse the
ribbons, while the discs are embedded.
Thus, without the ribbons, we get the $n$ component unlink.

Note that $\tau(T)$ is similar to $\phi(T)$ with the added stitching
$(t_1,~t_{2n})$.
The only difference between the two is the ribbons in the latter that are not
present in the former (see \figFinal).
However, we point out that removing the ribbons from a ribbon presentation
leaves just the components, which are equivalent to an unlink.
Since $\tau(R')$ is equivalent to $K$ with the ribbons removed, $\tau(R')$ is an
unlink, so $(\beta\circ\tau)(T)$ is an unlink as well.
Therefore, $T$ satisfies the conditions of the theorem.
\end{proof}

\papersec{Symmetric Unions}

The proof of \thmRibbon gives us a way of finding a tangle satisfying the stated
properties.
All that it requires is a ribbon presentation of a knot.
However, we notice that it was easier to find such tangles for the first 21
ribbon knots by using their symmetric unions.

It is known that all 21 ribbon knots with 10 crossings or fewer can be drawn as
symmetric unions (see \cite{many} and \cite{one}).
This helps us find tangles satisfying the conditions of \thmRibbon.
This is due to the fact that the easiest way to construct a tangle satisfying
these conditions is to cut the knot in several places, which will almost
automatically satisfy the condition that the full closure of the tangle must be
equivalent to the knot.
The other condition, that a particular closure of the tangle must create the
unlink, is harder to satisfy.
Note that the fusing number of each of the first 21 ribbon knots is 1 and the
knots' tangles can be drawn with only 4 strands.

These tangles are obtained from the original knots by cutting the knot in four
places.
The stitchings that would undo our cuts are $(t_1,~t_4)$, $(b_1,~b_2)$,
$(t_2,~t_3)$, and $(b_3,~b_4)$.
From this alone, it is clear that the resulting tangle will satisfy the second
condition of \thmRibbon, that its full closure with one additional stitching
will be the original knot.
Now we just need to ensure that the top and bottom closures of our tangle
result in the unlink.

The locations of the cuts that result in the ends $b_1$, $b_2$, $b_3$, and
$b_4$ are not important, but the two cuts resulting in the top ends are.
In a symmetric union, the central axis contains crossings and at least two
bridges.
One of these bridges is moved up until it is the uppermost part of the knot on
the axis.\newsavebox{\knotR}\sbox{\knotR}{\reflectbox{\size{9}{$R$}}}\\

\paperfig{Question}{\def\svgwidth{0.4\columnwidth}
{\centering\input{questionone.pdf_tex}}\\

\vspace{-6em}\hspace{0.3in}\hspace{14.8ex}\size{20}{$\equiv$}\\

\vspace{-6em}\hfill\def\svgwidth{0.4\columnwidth}
{\centering\input{questiontwo.pdf_tex}}\\

\hspace{0.3in}\hspace{15.1ex}\size{20}{$\Downarrow$}\\

\def\svgwidth{0.4\columnwidth}
{\centering\input{questionthree.pdf_tex}}\\

\vspace{-6em}\hspace{0.3in}\hspace{14.8ex}\size{20}{$\equiv$}\\

\vspace{-6.7em}\hfill\def\svgwidth{0.4\columnwidth}
{\centering\input{questionfour.pdf_tex}}}
{A visual representation of an implication involving tangles.
This shows that if a knot is a symmetric union represented by a tangle $R$, its
reflection \usebox{\knotR}, a series of crossings and bridges $X$, and the
appropriate connections, then cutting and reconnecting the strands in a certain
manner results in the closed untangle $U$.
Here, the knot has been cut just outside of the upper bridge and the ends were
connected symmetrically across the axis.
Afterwards, the inner connection was pulled under the original bridge.}\\

After moving the bridge to the top of the knot diagram, we follow the bridge
both ways until we reach a strand that crosses the strand of the bridge.
Then, we move along those strands to the outside of the knot and make our two
cuts there (see \figQuestion).
These cuts form the four upper ends of the resulting tangle for any of our 21
ribbon knots.
This raises the question of whether this can be generalized to all ribbon
knots.

\paperqtn{Cuts}
{Could the locations of the cuts between $(t_1,~t_4)$ and $(t_2,~t_3)$ for any
ribbon knot in a symmetric union presentation be generalized in order to be
along the two strands that are first to cross a bridge?}

It is known that all symmetric union presentations are ribbon knots, but it is
unknown whether every ribbon knot has a symmetric union presentation.
The answers we seek might depend on whether or not this is true.

To show how the tangle is obtained from a symmetric union, we demonstrate the
process for the simplest ribbon knot, $6_1$.
We start by drawing $6_1$ as a symmetric union.

\papersvg{ConvertOne}{convertone}{The ribbon knot $6_1$ drawn as a symmetric
union. The axis of symmetry is vertical.}\\

We need to verify that we can obtain the unlink by making the cuts and
stitchings between the four top ends.
As mentioned, for the 21 ribbon knots with 10 crossings or fewer we can always
make our cuts in a generalized location (see \qtnCuts).
Here, this location is on the two handles in the upper right and upper left
corners.
Having cut there, we stitch the ends that will later become $(t_1,~t_2)$ and
$(t_3,~t_4)$.\\

\papersvg{ConvertTwo}{converttwo}{The knot $6_1$ with the connections
$(t_1,~t_4)$ and $(t_2,~t_3)$ severed and stitchings between $(t_1,~t_2)$ and
$(t_3,~t_4)$.}\\

We now start to untangle this knot by untwisting the crossing in the middle of
the knot.

\papersvg{ConvertThree}{convertthree}{The restitched knot $6_1$ almost
completely untangled.}\\

Now, it is clear that the knot completely untwists into the 2-component unlink,
which is exactly what we wanted.

\papersvg{ConvertFour}{convertfour}{The 2-component unlink.}\\

Now that we know that our cuts were made in the right place, we can label
the ends both along the top and bottom with the labels 1 through 4.
This will let us know where to move the ends to construct our tangle.

\papersvg{ConvertedOne}{convertedone}{The symmetric union of $6_1$ with the ends
both along the top and bottom labeled 1 through 4.}

By pulling on the strands of the knot and adding a few kinks and crossings, we
arrange the ends along the top and bottom in the proper order, creating the
tangle we seek.\\

\papersvg{ConvertedTwo}{convertedtwo}{The constructed tangle for the knot
$6_1$.}\\

\papersec{Database}

By applying the process described in \thmRibbon to each of the 21 ribbon knots,
we were able to find the corresponding tangles.
We write down the tangle information of each of them by enumerating all of the
edges from a common starting point.
We use $t_1$ and then proceed downwards.
The edge descending from $t_1$ is labeled 1, and the value increments each time
it passes through a crossing.
We then write down the planar diagram information for the knot.\\

\paperfig{Crossing}
{\hspace{0.3in}\hspace{0.275\columnwidth}$k$\hspace{0.2\columnwidth}$j$
\vspace{-0.5em}\begin{center}\def\svgwidth{0.2\columnwidth}
{\centering\input{crossing.pdf_tex}}\end{center}
\vspace{-0.5em}
\hspace{0.3in}\hspace{0.275\columnwidth}$l$\hspace{0.2\columnwidth}$i$

\hspace{0.3in}\hspace{0.31\columnwidth}$X_{i,~j,~k,~l}$}
{A right-handed crossing labeled in planar diagram notation.
The lower incoming edge is labeled $i$ and then the remaining three are
labeled $j$, $k$, and $l$, proceeding counterclockwise from $i$.
The crossing is labeled as $X_{i,~j,~k,~l}$.}\\

Each crossing is labeled with $X_{i,~j,~k,~l}$ where $i$ is the index of the
lower incoming edge and the remaining indices proceed counterclockwise (see
\figCrossing).
The flow of the tangle is down the first strand, up the second, down the third,
and then up the fourth.\\

\papersvg{Example}{example}
{This depicts all the segments of knot $6_1$ as numbered.
The numbering starts from $t_1$ and proceeds downwards.
The tangle values of $6_1$, the values of the edges along the three connections
from strands 1 to 4, are (9, 13, 17).}\\

We also need to know where the strands of the tangle have been cut, so we need
to know the three numbers given to the connections $(b_1,~b_2),~(t_2,~t_3)$, and
$(b_3,~b_4)$ in that order.
We do not need to specify the value for $(t_1,~t_4)$ as it is always 1.
For the knot $6_1$, the tangle values are (9, 13, 17) (see \figExample).

We present these tangles for all 21 ribbon knots with 10 crossings or fewer,
drawn with the full closure and the connection $(t_1,~t_4)$ added, effectively
recreating the original knots.
We also present the planar diagram notation and tangle values of each of the
knots as they appear in the table (see below).
All of the data in the second table is also available online at\\
\url{http://tiny.cc/RibbonKnotDB}.

\papersec{Rectifying a Knot Table}

A symmetric union of a knot is not the only way to guarantee that we can find
where to cut the knot to create our tangles.
It is easy to show that this can always be done for any ribbon knot that has a
2-component ribbon presentation.
To do this, we cut each component in one place, as well as cut right across the
width of the ribbon, thereby making two cuts there and four cuts in total.

The author was originally planning on using \cite{knots} to accomplish this, as
it includes a table of 2-component ribbon presentations of all of the ribbon
knots with 10 crossings or fewer.
However, since both the symmetric unions in \cite{one} and \cite{many} and the
ribbon presentations in \cite{knots} contained far more crossings than the
minimal crossing number of the respective knots, they were verified for accuracy
by computing both their Alexander and Jones polynomials.

It was found that the Alexander polynomials of all of the knots in \cite{knots}
were correct, but the Jones polynomials of four of the knots were not.
These knots are $9_{27}$, $10_{42}$, $10_{75}$, and $10_{99}$.
The reason for this might be that a double delta move was accidentally performed
while in search of a ribbon presentation as the double delta move preserves the
Alexander polynomial.

\papersec{Acknowledgements}

The author is indebted to his advisor Professor Dror Bar-Natan for posing the problem as well as his tireless guidance and support during this project.
The author is also grateful to Akio Kawauchi for fruitful discussions on ribbon
knots and presentations.
Finally, the author thanks Christoph Lamm for his remarks and suggestions about
the questions and examples in this paper.

The research was partially supported by the National Sciences and Engineering
Research Council of Canada.

\papersec{References}

\end{paper}

\setlength{\tabcolsep}{12pt}
\begin{tabular}{cccc}
\def\svgwidth{0.17\columnwidth}
{\centering\input{6_1.pdf_tex}}&\def\svgwidth{0.17\columnwidth}
{\centering\input{8_8.pdf_tex}}&
\def\svgwidth{0.17\columnwidth}
{\centering\input{8_9.pdf_tex}}&\def\svgwidth{0.17\columnwidth}
{\centering\input{8_20.pdf_tex}}\\
$6_1$&$8_8$&$8_9$&$8_{20}$\\
&&&\\
\def\svgwidth{0.17\columnwidth}
{\centering\input{9_27.pdf_tex}}&\def\svgwidth{0.17\columnwidth}
{\centering\input{9_41.pdf_tex}}&
\def\svgwidth{0.17\columnwidth}
{\centering\input{9_46.pdf_tex}}&\def\svgwidth{0.17\columnwidth}
{\centering\input{10_3.pdf_tex}}\\
$9_{27}$&$9_{41}$&$9_{46}$&$10_3$\\
&&&\\
\def\svgwidth{0.17\columnwidth}
{\centering\input{10_22.pdf_tex}}&\def\svgwidth{0.17\columnwidth}
{\centering\input{10_35.pdf_tex}}&
\def\svgwidth{0.17\columnwidth}
{\centering\input{10_42.pdf_tex}}&\def\svgwidth{0.17\columnwidth}
{\centering\input{10_48.pdf_tex}}\\
$10_{22}$&$10_{35}$&$10_{42}$&$10_{48}$\\
&&&\\
\def\svgwidth{0.17\columnwidth}
{\centering\input{10_75.pdf_tex}}&\def\svgwidth{0.17\columnwidth}
{\centering\input{10_87.pdf_tex}}&
\def\svgwidth{0.17\columnwidth}
{\centering\input{10_99.pdf_tex}}&\def\svgwidth{0.17\columnwidth}
{\centering\input{10_123.pdf_tex}}\\
$10_{75}$&$10_{87}$&$10_{99}$&$10_{123}$\\
&&&\\
\def\svgwidth{0.17\columnwidth}
{\centering\input{10_129.pdf_tex}}&\def\svgwidth{0.17\columnwidth}
{\centering\input{10_137.pdf_tex}}&
\def\svgwidth{0.17\columnwidth}
{\centering\input{10_140.pdf_tex}}&\def\svgwidth{0.17\columnwidth}
{\centering\input{10_153.pdf_tex}}\\
$10_{129}$&$10_{137}$&$10_{140}$&$10_{153}$\\
&&&\\
\def\svgwidth{0.17\columnwidth}
{\centering\input{10_155.pdf_tex}}&&&\\
$10_{155}$&&&
\end{tabular}

\pagebreak

\begin{table}[h]
\begin{tabularx}{\textwidth}{|c|X|c|}\hline
\brokensize{9}{Knot}&\brokensize{9}{Planar Diagram Notation}&
\brokensize{9}{Tangle Values}\\\hline
\brokensize{9}{$6_1$}&\brokensize{9}{$X_{2,~8,~3,~7}$ $X_{3,~10,~4,~11}$
$X_{5,~14,~6,~15}$ $X_{8,~20,~9,~19}$ $X_{11,~6,~12,~7}$ $X_{15,~4,~16,~5}$
$X_{17,~16,~18,~17}$ $X_{18,~10,~19,~9}$ $X_{20,~2,~21,~1}$ $X_{21,~12,~22,~13}$
$X_{22,~14,~1,~13}$}&\brokensize{9}{(9,~13,~17)}\\\hline
\brokensize{9}{$8_8$}&\brokensize{9}{$X_{2,~10,~3,~9}$ $X_{4,~24,~5,~23}$
$X_{6,~12,~7,~11}$ $X_{7,~14,~8,~15}$ $X_{10,~4,~11,~3}$ $X_{13,~18,~14,~19}$
$X_{15,~8,~16,~9}$ $X_{19,~12,~20,~13}$ $X_{21,~20,~22,~21}$ $X_{22,~6,~23,~5}$
$X_{24,~2,~25,~1}$ $X_{25,~16,~26,~17}$
$X_{26,~18,~1,~17}$}&\brokensize{9}{(5,~17,~21)}\\\hline
\brokensize{9}{$8_9$}&\brokensize{9}{$X_{2,~13,~3,~14}$ $X_{4,~18,~5,~17}$
$X_{9,~8,~10,~9}$ $X_{11,~30,~12,~31}$ $X_{16,~27,~17,~28}$ $X_{18,~30,~19,~29}$
$X_{19,~23,~20,~22}$ $X_{20,~8,~21,~7}$ $X_{21,~6,~22,~7}$ $X_{23,~11,~24,~10}$
$X_{25,~12,~26,~13}$ $X_{26,~4,~27,~3}$ $X_{28,~5,~29,~6}$ $X_{31,~25,~32,~24}$
$X_{32,~2,~33,~1}$ $X_{33,~14,~34,~15}$
$X_{34,~16,~1,~15}$}&\brokensize{9}{(9,~15,~21)}\\\hline
\brokensize{9}{$8_{20}$}&\brokensize{9}{$X_{2,~9,~3,~10}$ $X_{3,~14,~4,~15}$
$X_{7,~6,~8,~7}$ $X_{8,~23,~9,~24}$ $X_{11,~19,~12,~18}$ $X_{13,~4,~14,~5}$
$X_{15,~11,~16,~10}$ $X_{19,~13,~20,~12}$ $X_{21,~20,~22,~21}$
$X_{22,~5,~23,~6}$ $X_{24,~2,~25,~1}$ $X_{25,~16,~26,~17}$
$X_{26,~18,~1,~17}$}&\brokensize{9}{(7,~17,~21)}\\\hline
\brokensize{9}{$9_{27}$}&\brokensize{9}{$X_{3,~27,~4,~26}$ $X_{6,~23,~7,~24}$
$X_{8,~16,~9,~15}$ $X_{9,~2,~10,~3}$ $X_{11,~23,~12,~22}$ $X_{12,~6,~13,~5}$
$X_{14,~19,~15,~20}$ $X_{18,~8,~19,~7}$ $X_{21,~4,~22,~5}$ $X_{24,~14,~25,~13}$
$X_{25,~20,~26,~21}$ $X_{27,~10,~28,~11}$ $X_{28,~2,~29,~1}$
$X_{29,~16,~30,~17}$ $X_{30,~18,~1,~17}$}&\brokensize{9}{(11,~17,~23)}\\\hline
\brokensize{9}{$9_{41}$}&\brokensize{9}{$X_{3,~18,~4,~19}$ $X_{4,~24,~5,~23}$
$X_{7,~14,~8,~15}$ $X_{8,~27,~9,~28}$ $X_{11,~21,~12,~20}$ $X_{12,~1,~13,~2}$
$X_{13,~27,~14,~26}$ $X_{15,~6,~16,~7}$ $X_{17,~30,~18,~31}$ $X_{22,~10,~23,~9}$
$X_{24,~30,~25,~29}$ $X_{25,~17,~26,~16}$ $X_{28,~5,~29,~6}$ $X_{31,~3,~32,~2}$
$X_{32,~19,~33,~20}$ $X_{33,~10,~34,~11}$
$X_{34,~22,~1,~21}$}&\brokensize{9}{(13,~21,~27)}\\\hline
\brokensize{9}{$9_{46}$}&\brokensize{9}{$X_{2,~10,~3,~9}$ $X_{3,~14,~4,~15}$
$X_{5,~12,~6,~13}$ $X_{7,~18,~8,~19}$ $X_{10,~24,~11,~23}$ $X_{13,~4,~14,~5}$
$X_{15,~8,~16,~9}$ $X_{19,~6,~20,~7}$ $X_{21,~20,~22,~21}$ $X_{22,~12,~23,~11}$
$X_{24,~2,~25,~1}$ $X_{25,~16,~26,~17}$
$X_{26,~18,~1,~17}$}&\brokensize{9}{(11,~17,~21)}\\\hline
\brokensize{9}{$10_3$}&\brokensize{9}{$X_{2,~13,~3,~14}$ $X_{3,~18,~4,~19}$
$X_{5,~16,~6,~17}$ $X_{8,~10,~9,~9}$ $X_{10,~24,~11,~23}$ $X_{12,~22,~13,~21}$
$X_{15,~6,~16,~7}$ $X_{17,~4,~18,~5}$ $X_{19,~15,~20,~14}$ $X_{20,~1,~21,~2}$
$X_{22,~12,~23,~11}$ $X_{24,~8,~1,~7}$}&\brokensize{9}{(9,~15,~1)}\\\hline
\brokensize{9}{$10_{22}$}&\brokensize{9}{$X_{3,~17,~4,~16}$ $X_{8,~23,~9,~24}$
$X_{12,~21,~13,~22}$ $X_{14,~8,~15,~7}$ $X_{17,~31,~18,~30}$ $X_{19,~3,~20,~2}$
$X_{20,~9,~21,~10}$ $X_{22,~13,~23,~14}$ $X_{24,~16,~25,~15}$ $X_{25,~6,~26,~7}$
$X_{27,~26,~28,~27}$ $X_{28,~6,~29,~5}$ $X_{29,~4,~30,~5}$ $X_{31,~19,~32,~18}$
$X_{32,~2,~33,~1}$ $X_{33,~10,~34,~11}$
$X_{34,~12,~1,~11}$}&\brokensize{9}{(5,~11,~27)}\\\hline
\brokensize{9}{$10_{35}$}&\brokensize{9}{$X_{3,~12,~4,~13}$ $X_{4,~32,~5,~31}$
$X_{6,~30,~7,~29}$ $X_{8,~14,~9,~13}$ $X_{10,~20,~11,~19}$ $X_{11,~2,~12,~3}$
$X_{15,~24,~16,~25}$ $X_{17,~22,~18,~23}$ $X_{18,~10,~19,~9}$
$X_{23,~16,~24,~17}$ $X_{25,~14,~26,~15}$ $X_{27,~26,~28,~27}$
$X_{28,~8,~29,~7}$ $X_{30,~6,~31,~5}$ $X_{32,~2,~33,~1}$ $X_{33,~20,~34,~21}$
$X_{34,~22,~1,~21}$}&\brokensize{9}{(7,~21,~27)}\\\hline
\brokensize{9}{$10_{42}$}&\brokensize{9}{$X_{3,~39,~4,~38}$ $X_{5,~14,~6,~15}$
$X_{6,~35,~7,~36}$ $X_{9,~13,~10,~12}$ $X_{13,~25,~14,~24}$ $X_{16,~21,~17,~22}$
$X_{20,~33,~21,~34}$ $X_{22,~38,~23,~37}$ $X_{23,~4,~24,~5}$ $X_{25,~9,~26,~8}$
$X_{27,~26,~28,~27}$ $X_{28,~12,~29,~11}$ $X_{29,~10,~30,~11}$
$X_{31,~3,~32,~2}$ $X_{32,~17,~33,~18}$ $X_{34,~7,~35,~8}$ $X_{36,~16,~37,~15}$
$X_{39,~31,~40,~30}$ $X_{40,~2,~41,~1}$ $X_{41,~18,~42,~19}$
$X_{42,~20,~1,~19}$}&\brokensize{9}{(11,~19,~27)}\\\hline
\brokensize{9}{$10_{48}$}&\brokensize{9}{$X_{4,~14,~5,~13}$ $X_{5,~27,~6,~26}$
$X_{7,~29,~8,~28}$ $X_{9,~14,~10,~15}$ $X_{11,~22,~12,~23}$ $X_{12,~4,~13,~3}$
$X_{17,~16,~18,~17}$ $X_{18,~33,~19,~34}$ $X_{20,~35,~21,~36}$
$X_{21,~10,~22,~11}$ $X_{23,~3,~24,~2}$ $X_{27,~7,~28,~6}$ $X_{29,~9,~30,~8}$
$X_{31,~30,~32,~31}$ $X_{32,~15,~33,~16}$ $X_{34,~19,~35,~20}$
$X_{36,~2,~37,~1}$ $X_{37,~24,~38,~25}$
$X_{38,~26,~1,~25}$}&\brokensize{9}{(17,~25,~31)}\\\hline
\brokensize{9}{$10_{75}$}&\brokensize{9}{$X_{4,~14,~5,~13}$ $X_{5,~27,~6,~26}$
$X_{7,~29,~8,~28}$ $X_{9,~14,~10,~15}$ $X_{11,~22,~12,~23}$ $X_{12,~4,~13,~3}$
$X_{17,~16,~18,~17}$ $X_{18,~33,~19,~34}$ $X_{20,~35,~21,~36}$
$X_{21,~10,~22,~11}$ $X_{23,~3,~24,~2}$ $X_{27,~7,~28,~6}$ $X_{29,~9,~30,~8}$
$X_{31,~30,~32,~31}$ $X_{32,~15,~33,~16}$ $X_{34,~19,~35,~20}$
$X_{36,~2,~37,~1}$ $X_{37,~24,~38,~25}$
$X_{38,~26,~1,~25}$}&\brokensize{9}{(9,~17,~25)}\\\hline
\brokensize{9}{$10_{87}$}&\brokensize{9}{$X_{4,~22,~5,~21}$ $X_{6,~14,~7,~13}$
$X_{7,~24,~8,~25}$ $X_{12,~36,~13,~35}$ $X_{15,~38,~16,~39}$ $X_{16,~4,~17,~3}$
$X_{20,~33,~21,~34}$ $X_{23,~15,~24,~14}$ $X_{25,~10,~26,~11}$
$X_{27,~26,~28,~27}$ $X_{28,~10,~29,~9}$ $X_{29,~8,~30,~9}$ $X_{31,~3,~32,~2}$
$X_{32,~17,~33,~18}$ $X_{34,~12,~35,~11}$ $X_{36,~5,~37,~6}$
$X_{37,~22,~38,~23}$ $X_{39,~31,~40,~30}$ $X_{40,~2,~41,~1}$
$X_{41,~18,~42,~19}$ $X_{42,~20,~1,~19}$}&\brokensize{9}{(9,~19,~27)}\\\hline
\brokensize{9}{$10_{99}$}&\brokensize{9}{$X_{2,~29,~3,~30}$ $X_{4,~13,~5,~14}$
$X_{7,~36,~8,~37}$ $X_{11,~21,~12,~20}$ $X_{12,~3,~13,~4}$ $X_{14,~9,~15,~10}$
$X_{16,~26,~17,~25}$ $X_{19,~11,~20,~10}$ $X_{21,~31,~22,~30}$
$X_{24,~18,~25,~17}$ $X_{27,~9,~28,~8}$ $X_{28,~5,~29,~6}$ $X_{31,~19,~32,~18}$
$X_{32,~15,~33,~16}$ $X_{33,~27,~34,~26}$ $X_{35,~34,~36,~35}$
$X_{37,~6,~38,~7}$ $X_{38,~2,~39,~1}$ $X_{39,~22,~40,~23}$
$X_{40,~24,~1,~23}$}&\brokensize{9}{(7,~23,~35)}\\\hline
\brokensize{9}{$10_{123}$}&\brokensize{9}{$X_{2,~23,~3,~24}$ $X_{5,~14,~6,~15}$
$X_{7,~17,~8,~16}$ $X_{8,~3,~9,~4}$ $X_{10,~34,~11,~33}$ $X_{15,~26,~16,~27}$
$X_{17,~25,~18,~24}$ $X_{20,~14,~21,~13}$ $X_{22,~9,~23,~10}$ $X_{25,~7,~26,~6}$
$X_{27,~4,~28,~5}$ $X_{28,~22,~29,~21}$ $X_{29,~12,~30,~13}$
$X_{31,~30,~32,~31}$ $X_{32,~12,~33,~11}$ $X_{34,~2,~35,~1}$
$X_{35,~18,~36,~19}$ $X_{36,~20,~1,~19}$}&\brokensize{9}{(11,~19,~31)}\\\hline
\brokensize{9}{$10_{129}$}&\brokensize{9}{$X_{2,~19,~3,~20}$ $X_{4,~12,~5,~11}$
$X_{6,~28,~7,~27}$ $X_{8,~14,~9,~13}$ $X_{9,~16,~10,~17}$ $X_{12,~6,~13,~5}$
$X_{15,~22,~16,~23}$ $X_{17,~10,~18,~11}$ $X_{18,~3,~19,~4}$
$X_{23,~14,~24,~15}$ $X_{25,~24,~26,~25}$ $X_{26,~8,~27,~7}$ $X_{28,~2,~29,~1}$
$X_{29,~20,~30,~21}$ $X_{30,~22,~1,~21}$}&\brokensize{9}{(7,~21,~25)}\\\hline
\brokensize{9}{$10_{137}$}&\brokensize{9}{$X_{3,~27,~4,~26}$ $X_{8,~13,~9,~14}$
$X_{12,~18,~13,~17}$ $X_{15,~24,~16,~25}$ $X_{16,~7,~17,~8}$ $X_{18,~10,~19,~9}$
$X_{19,~2,~20,~3}$ $X_{21,~5,~22,~4}$ $X_{22,~5,~23,~6}$ $X_{23,~7,~24,~6}$
$X_{25,~14,~26,~15}$ $X_{27,~20,~28,~21}$ $X_{28,~2,~29,~1}$
$X_{29,~10,~30,~11}$ $X_{30,~12,~1,~11}$}&\brokensize{9}{(5,~11,~23)}\\\hline
\brokensize{9}{$10_{140}$}&\brokensize{9}{$X_{2,~11,~3,~12}$ $X_{3,~18,~4,~19}$
$X_{5,~16,~6,~17}$ $X_{9,~8,~10,~9}$ $X_{10,~27,~11,~28}$ $X_{13,~23,~14,~22}$
$X_{15,~6,~16,~7}$ $X_{17,~4,~18,~5}$ $X_{19,~13,~20,~12}$ $X_{23,~15,~24,~14}$
$X_{25,~24,~26,~25}$ $X_{26,~7,~27,~8}$ $X_{28,~2,~29,~1}$ $X_{29,~20,~30,~21}$
$X_{30,~22,~1,~21}$}&\brokensize{9}{(9,~21,~25)}\\\hline
\brokensize{9}{$10_{153}$}&\brokensize{9}{$X_{3,~27,~4,~26}$ $X_{4,~22,~5,~21}$
$X_{7,~6,~8,~7}$ $X_{9,~3,~10,~2}$ $X_{10,~15,~11,~16}$ $X_{13,~24,~14,~25}$
$X_{14,~19,~15,~20}$ $X_{18,~11,~19,~12}$ $X_{22,~6,~23,~5}$
$X_{23,~12,~24,~13}$ $X_{25,~20,~26,~21}$ $X_{27,~9,~28,~8}$ $X_{28,~2,~29,~1}$
$X_{29,~16,~30,~17}$ $X_{30,~18,~1,~17}$}&\brokensize{9}{(7,~17,~23)}\\\hline
\brokensize{9}{$10_{155}$}&\brokensize{9}{$X_{4,~18,~5,~17}$ $X_{9,~8,~10,~9}$
$X_{11,~30,~12,~31}$ $X_{14,~2,~15,~1}$ $X_{16,~27,~17,~28}$
$X_{18,~30,~19,~29}$ $X_{20,~8,~21,~7}$ $X_{21,~6,~22,~7}$ $X_{22,~19,~23,~20}$
$X_{23,~11,~24,~10}$ $X_{25,~12,~26,~13}$ $X_{26,~4,~27,~3}$ $X_{28,~5,~29,~6}$
$X_{31,~25,~32,~24}$ $X_{32,~13,~33,~14}$ $X_{33,~3,~34,~2}$
$X_{34,~16,~1,~15}$}&\brokensize{9}{(9,~15,~21)}\\\hline
\end{tabularx}
\end{table}
\end{document}

%% file: header.tex
\usepackage{amsmath}
\usepackage{amssymb}
\usepackage{amsthm}
\usepackage{capt-of}
\usepackage{caption}
\usepackage[strict]{changepage}
\usepackage{chngcntr}
\usepackage[americanvoltage,siunitx]{circuitikz}
\usepackage{color,colortbl}
\usepackage{etoolbox}
\usepackage{fancyhdr}
\usepackage[T1]{fontenc}
\usepackage{gensymb}
\usepackage[margin=1in]{geometry}
\usepackage{graphicx}
\usepackage{hyperref}
\usepackage{import}
\usepackage{indentfirst}
\usepackage{mathptmx}
\usepackage{mathrsfs}
\usepackage{multicol}
\usepackage{multirow}
\usepackage{needspace}
\usepackage{pgfplots}
\usepackage{pgfplotstable}
\usepackage{setspace}
\usepackage{siunitx}
\usepackage{tabu}
\usepackage{tabularx}
\usepackage{tikz}
\usepackage{xspace}

\patchcmd{\thebibliography}{\section*{\refname}}{\vspace{-1em}}{}{}

\usepgfplotslibrary{external}
\usetikzlibrary{positioning,matrix,shapes,chains,arrows}
\newcommand\svgc[1]{\def\svgwidth{\columnwidth}
{\centering\input{#1.pdf_tex}}}

\newcommand\size[2]{{\fontsize{#1pt}{#1pt}\selectfont#2}}
\newcommand\brokensize[2]{\fontsize{#1pt}{#1pt}\selectfont#2}

\newcounter{paperthmamount}
\newcommand\theorems[0]{
\theoremstyle{remark}
\newtheorem{claim}[subsection]{Claim}
\theoremstyle{plain}
\newtheorem{conjecture}[subsection]{Conjecture}
\theoremstyle{plain}
\newtheorem{corollary}[subsection]{Corollary}
\theoremstyle{definition}
\newtheorem{definition}[subsection]{Definition}
\theoremstyle{plain}
\newtheorem{lemma}[subsection]{Lemma}
\theoremstyle{remark}
\newtheorem{proposition}[subsection]{Proposition}
\theoremstyle{remark}
\newtheorem{remark}[subsection]{Remark}
\theoremstyle{plain}
\newtheorem{theorem}[subsection]{Theorem}
\theoremstyle{definition}
\newtheorem{question}[subsection]{Question}
\newcommand\paperclm[2]
{\begin{claim}\global\expandafter\edef
\csname clm##1\endcsname{Claim \thesubsection\noexpand\xspace}
##2\end{claim}}
\newcommand\papercnj[2]
{\begin{conjecture}\global\expandafter\edef
\csname cnj##1\endcsname{Conjecture \thesubsection\noexpand\xspace}
##2\end{conjecture}}
\newcommand\papercor[2]
{\begin{corollary}\global\expandafter\edef
\csname cor##1\endcsname{Corollary \thesubsection\noexpand\xspace}
##2\end{corollary}}
\newcommand\paperdef[2]
{\begin{definition}\global\expandafter\edef
\csname def##1\endcsname{Definition \thesubsection\noexpand\xspace}
##2\end{definition}}
\newcommand\paperlem[2]
{\begin{lemma}\global\expandafter\edef
\csname lem##1\endcsname{Lemma \thesubsection\noexpand\xspace}
##2\end{lemma}}
\newcommand\paperprp[2]
{\begin{proposition}\global\expandafter\edef
\csname prp##1\endcsname{Proposition \thesubsection\noexpand\xspace}
##2\end{proposition}}
\newcommand\paperqtn[2]
{\begin{question}\global\expandafter\edef
\csname qtn##1\endcsname{Question \thesubsection\noexpand\xspace}
##2\end{question}}
\newcommand\paperrem[2]
{\begin{remark}\global\expandafter\edef
\csname rem##1\endcsname{Remark \thesubsection\noexpand\xspace}
##2\end{remark}}
\newcommand\paperthm[2]
{\begin{theorem}\global\expandafter\edef
\csname thm##1\endcsname{Theorem \thesubsection\noexpand\xspace}
##2\end{theorem}}}
\newcommand\subtheorems[0]{\stepcounter{paperthmamount}
\theoremstyle{remark}
\newtheorem{claim}[subsubsection]{Claim}
\theoremstyle{plain}
\newtheorem{conjecture}[subsubsection]{Conjecture}
\theoremstyle{plain}
\newtheorem{corollary}[subsubsection]{Corollary}
\theoremstyle{definition}
\newtheorem{definition}[subsubsection]{Definition}
\theoremstyle{plain}
\newtheorem{lemma}[subsubsection]{Lemma}
\theoremstyle{remark}
\newtheorem{proposition}[subsubsection]{Proposition}
\theoremstyle{remark}
\newtheorem{remark}[subsubsection]{Remark}
\theoremstyle{plain}
\newtheorem{theorem}[subsubsection]{Theorem}
\theoremstyle{definition}
\newtheorem{question}[subsubsection]{Question}
\newcommand\paperclm[2]
{\begin{claim}\global\expandafter\edef
\csname clm##1\endcsname{Claim \thesubsubsection\noexpand\xspace}
##2\end{claim}}
\newcommand\papercnj[2]
{\begin{conjecture}\global\expandafter\edef
\csname cnj##1\endcsname{Conjecture \thesubsubsection\noexpand\xspace}
##2\end{conjecture}}
\newcommand\papercor[2]
{\begin{corollary}\global\expandafter\edef
\csname cor##1\endcsname{Corollary \thesubsubsection\noexpand\xspace}
##2\end{corollary}}
\newcommand\paperdef[2]
{\begin{definition}\global\expandafter\edef
\csname def##1\endcsname{Definition \thesubsubsection\noexpand\xspace}
##2\end{definition}}
\newcommand\paperlem[2]
{\begin{lemma}\global\expandafter\edef
\csname lem##1\endcsname{Lemma \thesubsubsection\noexpand\xspace}
##2\end{lemma}}
\newcommand\paperprp[2]
{\begin{proposition}\global\expandafter\edef
\csname prp##1\endcsname{Proposition \thesubsubsection\noexpand\xspace}
##2\end{proposition}}
\newcommand\paperqtn[2]
{\begin{question}\global\expandafter\edef
\csname qtn##1\endcsname{Question \thesubsubsection\noexpand\xspace}
##2\end{question}}
\newcommand\paperrem[2]
{\begin{remark}\global\expandafter\edef
\csname rem##1\endcsname{Remark \thesubsubsection\noexpand\xspace}
##2\end{remark}}
\newcommand\paperthm[2]
{\begin{theorem}\global\expandafter\edef
\csname thm##1\endcsname{Theorem \thesubsubsection\noexpand\xspace}
##2\end{theorem}}}

\newcommand\papertitle[1]
{{\centering\fontsize{20pt}{20pt}\textsc{#1}\\\mbox{}\\}
\fancyhead[OC]{\fontsize{12pt}{12pt}\selectfont\textit{#1}}}
\newcounter{people}
\newcommand\paperauthtext[3]{{\centering\fontsize{12pt}{12pt}\selectfont
\textsc{#1}\\[-0.1em]{\fontsize{9pt}{9pt}\selectfont\textit{\ifx&#2&
\vspace{-1em}\else#2\fi}}\\\mbox{}\\
\fancyhead[EC]{\fontsize{12pt}{12pt}\selectfont\textit{#3}}}}
\newcommand\paperauth[2]{{\stepcounter{people}
\ifnum\value{people}=1
{\paperauthtext{#1}{#2}{#1}
\global\def\auth{#1\xspace}}
\else\ifnum\value{people}=2
{\paperauthtext{#1}{#2}{\auth and #1}}
\else{\paperauthtext{#1}{#2}{\auth et al}}\fi\fi}}
\newcommand\physics[0]{
\renewcommand\paperauthtext[4]{{\centering\fontsize{12pt}{12pt}\selectfont
\textsc{##1. ##2}\\[-0.1em]{\fontsize{9pt}{9pt}\selectfont\textit{\ifx&##3&
\vspace{-1em}\else##3\fi}}\\\mbox{}\\
\fancyhead[EC]{\fontsize{12pt}{12pt}\selectfont\textit{##4}}}}
\renewcommand\paperauth[3]{{\stepcounter{people}
\ifnum\value{people}=1
{\paperauthtext{##1}{##2}{##3}{##1. ##2}
\global\def\auth{##2\xspace}}
\else\ifnum\value{people}=2
{\paperauthtext{##1}{##2}{##3}{\auth and ##2}}
\else{\paperauthtext{##1}{##2}{##3}{\auth et al}}\fi\fi}}}

\makeatletter
\newenvironment{paperadjustwidth}[2]{
  \begin{list}{}{
    \setlength\partopsep\z@
    \setlength\topsep\z@
    \setlength\listparindent\parindent
    \setlength\parsep\parskip
    \linespread{0.75}\selectfont
    \@ifmtarg{#1}{\setlength{\leftmargin}{\z@}}
                 {\setlength{\leftmargin}{#1}}
    \@ifmtarg{#2}{\setlength{\rightmargin}{\z@}}
                 {\setlength{\rightmargin}{#2}}
    }
    \item[]}{\end{list}}
\makeatother

\newenvironment{paperabs}
{\begin{paperadjustwidth}{0.5in}{0.5in}\bgroup\fontsize{9pt}{9pt}\selectfont
\hspace{0.5in}}
{\egroup\end{paperadjustwidth}}

\newenvironment{paper}
{\begin{multicols*}{2}\bgroup\fontsize{12pt}{12pt}\selectfont}
{\egroup\end{multicols*}}

\newsavebox{\sourcebox}

\newcounter{paperseccounter}
\newcounter{papersubseccounter}[paperseccounter]
\newcommand\papersec[1]{\needspace{1in}
\stepcounter{paperseccounter}
\stepcounter{section}
\begin{center}\Roman{paperseccounter} \textsc{#1}\end{center}}

\newcounter{papereqcounter}
\newcommand\papereq[3]{{
\stepcounter{papereqcounter}
\mbox{}\vspace{-0.75em}
\begin{equation*}
#2
\tag*{\fontsize{12pt}{12pt}\selectfont
$\begin{array}{r}
\cr{\text{(\arabic{papereqcounter})}}
\cr{\fontsize{9pt}{9pt}\selectfont\textit{\ifx\\#3\\~\else(\fi#3\ifx\\#3\\~
\else)\fi}}
\end{array}$}
\end{equation*}
}
\expandafter\edef\csname eq#1\endcsname{(\arabic{papereqcounter})\noexpand
\xspace}}

\newcommand{\papervar}[3]
{&$#1$ & #2 \ifx\\#3\\~\else($\smash{\text{\si{\fi
#3\ifx\\#3\\~\else}}}$)\fi\\}
\newenvironment{paperwhere}
{\begin{minipage}{\columnwidth}
\bgroup\fontsize{9pt}{9pt}\selectfont Where:\vspace{2pt}\\\begin{tabular}
{rr@{ = }p{\linewidth}}}
{\end{tabular}\egroup\end{minipage}\vspace{5pt}}

\definecolor{LineGray}{gray}{0.5}
\newtabulinestyle{outer=2.25pt LineGray}
\newtabulinestyle{inner=0.75pt LineGray}
\newcolumntype{I}{X[-5,c]}
\newcolumntype{U}{@{}X[-5,r]@{$\pm$}X[-5,l]@{}}
\newcolumntype{C}{@{}X[-5,c]@{}}

\newcounter{papertableindexcounter}

\newcounter{paperfigurecounter}
\newcommand{\papercap}[2]{\bgroup\stepcounter{paperfigurecounter}
\captionof{figure}{\fontsize{9pt}{9pt}\selectfont
\hspace{0.3in}Fig.~\arabic{paperfigurecounter}.\quad#2}
\egroup\expandafter\edef
\csname fig#1\endcsname{Fig.~\arabic{paperfigurecounter}\noexpand\xspace}}

\newcommand\paperfig[3]{\noindent\begin{minipage}{\columnwidth}
#2\papercap{#1}{#3}\end{minipage}\expandafter\edef
\csname fig#1\endcsname{Fig.~\arabic{paperfigurecounter}\noexpand\xspace}}
\newcommand\papersvg[3]{\paperfig{#1}{\svgc{#2}}{#3}}

\newcommand{\comment}[1]{}

%% file: diagrams/clasp.pdf_tex
\begingroup%
  \makeatletter%
  \providecommand\color[2][]{%
    \errmessage{(Inkscape) Color is used for the text in Inkscape, but the package 'color.sty' is not loaded}%
    \renewcommand\color[2][]{}%
  }%
  \providecommand\transparent[1]{%
    \errmessage{(Inkscape) Transparency is used (non-zero) for the text in Inkscape, but the package 'transparent.sty' is not loaded}%
    \renewcommand\transparent[1]{}%
  }%
  \providecommand\rotatebox[2]{#2}%
  \ifx\svgwidth\undefined%
    \setlength{\unitlength}{596.47309171bp}%
    \ifx\svgscale\undefined%
      \relax%
    \else%
      \setlength{\unitlength}{\unitlength * \real{\svgscale}}%
    \fi%
  \else%
    \setlength{\unitlength}{\svgwidth}%
  \fi%
  \global\let\svgwidth\undefined%
  \global\let\svgscale\undefined%
  \makeatother%
  \begin{picture}(1,0.71913525)%
    \put(0,0){\includegraphics[width=\unitlength,page=1]{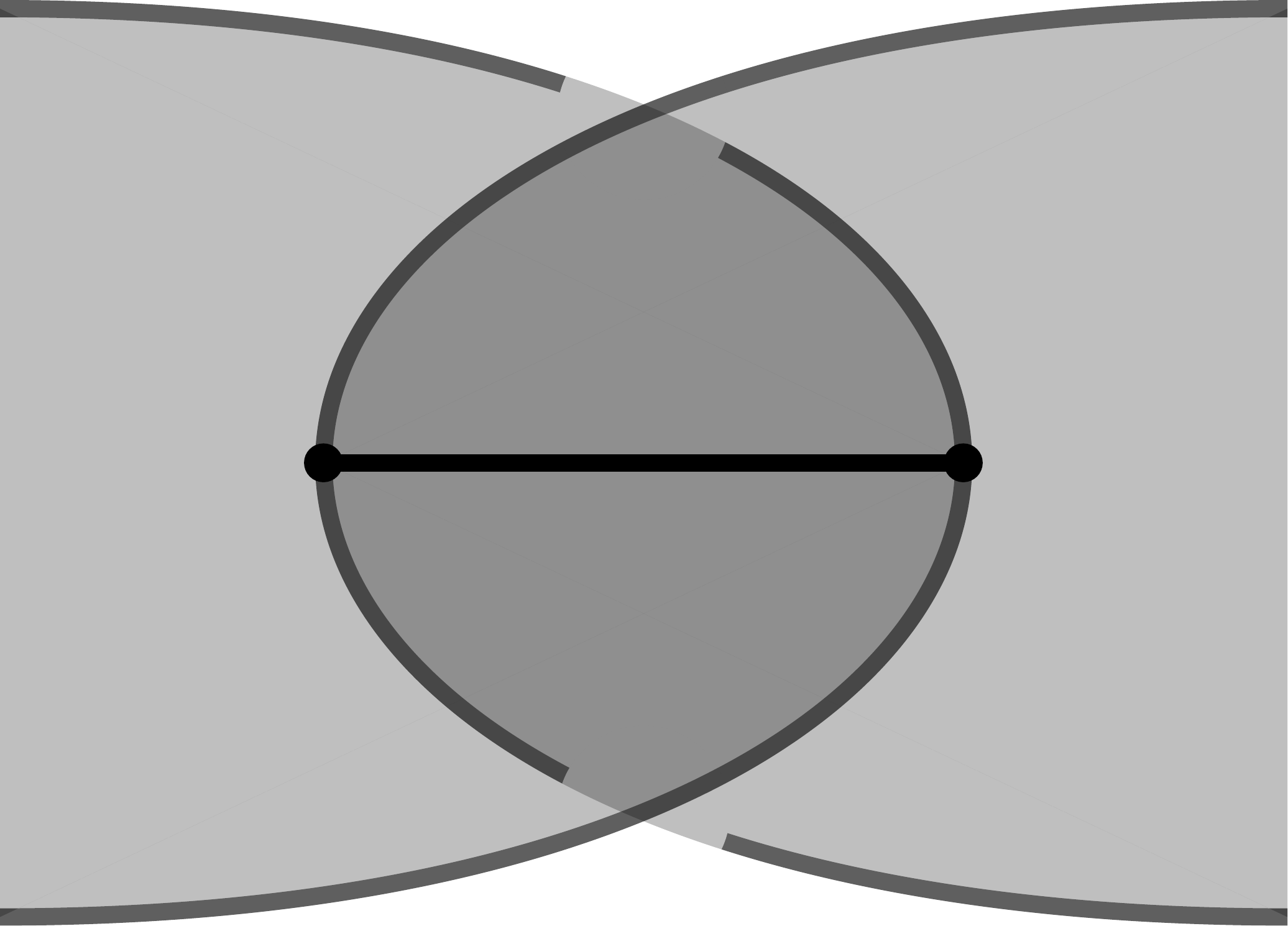}}%
  \end{picture}%
\endgroup%

%% file: diagrams/ribbon.pdf_tex
\begingroup%
  \makeatletter%
  \providecommand\color[2][]{%
    \errmessage{(Inkscape) Color is used for the text in Inkscape, but the package 'color.sty' is not loaded}%
    \renewcommand\color[2][]{}%
  }%
  \providecommand\transparent[1]{%
    \errmessage{(Inkscape) Transparency is used (non-zero) for the text in Inkscape, but the package 'transparent.sty' is not loaded}%
    \renewcommand\transparent[1]{}%
  }%
  \providecommand\rotatebox[2]{#2}%
  \ifx\svgwidth\undefined%
    \setlength{\unitlength}{602.08221354bp}%
    \ifx\svgscale\undefined%
      \relax%
    \else%
      \setlength{\unitlength}{\unitlength * \real{\svgscale}}%
    \fi%
  \else%
    \setlength{\unitlength}{\svgwidth}%
  \fi%
  \global\let\svgwidth\undefined%
  \global\let\svgscale\undefined%
  \makeatother%
  \begin{picture}(1,0.70998673)%
    \put(0,0){\includegraphics[width=\unitlength,page=1]{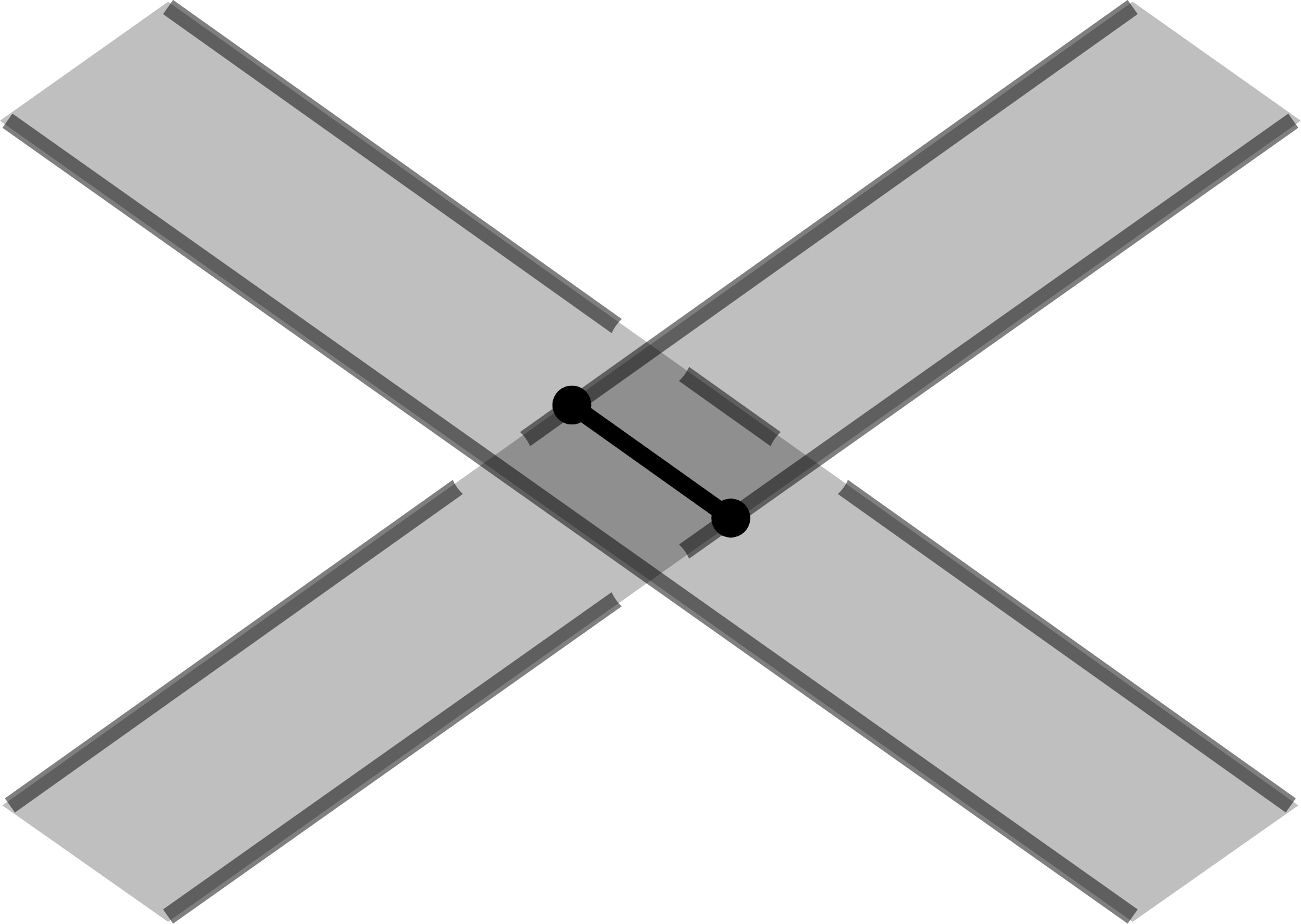}}%
  \end{picture}%
\endgroup%

%% file: diagrams/lemmaone.pdf_tex
\begingroup%
  \makeatletter%
  \providecommand\color[2][]{%
    \errmessage{(Inkscape) Color is used for the text in Inkscape, but the package 'color.sty' is not loaded}%
    \renewcommand\color[2][]{}%
  }%
  \providecommand\transparent[1]{%
    \errmessage{(Inkscape) Transparency is used (non-zero) for the text in Inkscape, but the package 'transparent.sty' is not loaded}%
    \renewcommand\transparent[1]{}%
  }%
  \providecommand\rotatebox[2]{#2}%
  \ifx\svgwidth\undefined%
    \setlength{\unitlength}{603.27558594bp}%
    \ifx\svgscale\undefined%
      \relax%
    \else%
      \setlength{\unitlength}{\unitlength * \real{\svgscale}}%
    \fi%
  \else%
    \setlength{\unitlength}{\svgwidth}%
  \fi%
  \global\let\svgwidth\undefined%
  \global\let\svgscale\undefined%
  \makeatother%
  \begin{picture}(1,0.78330714)%
    \put(0,0){\includegraphics[width=\unitlength,page=1]{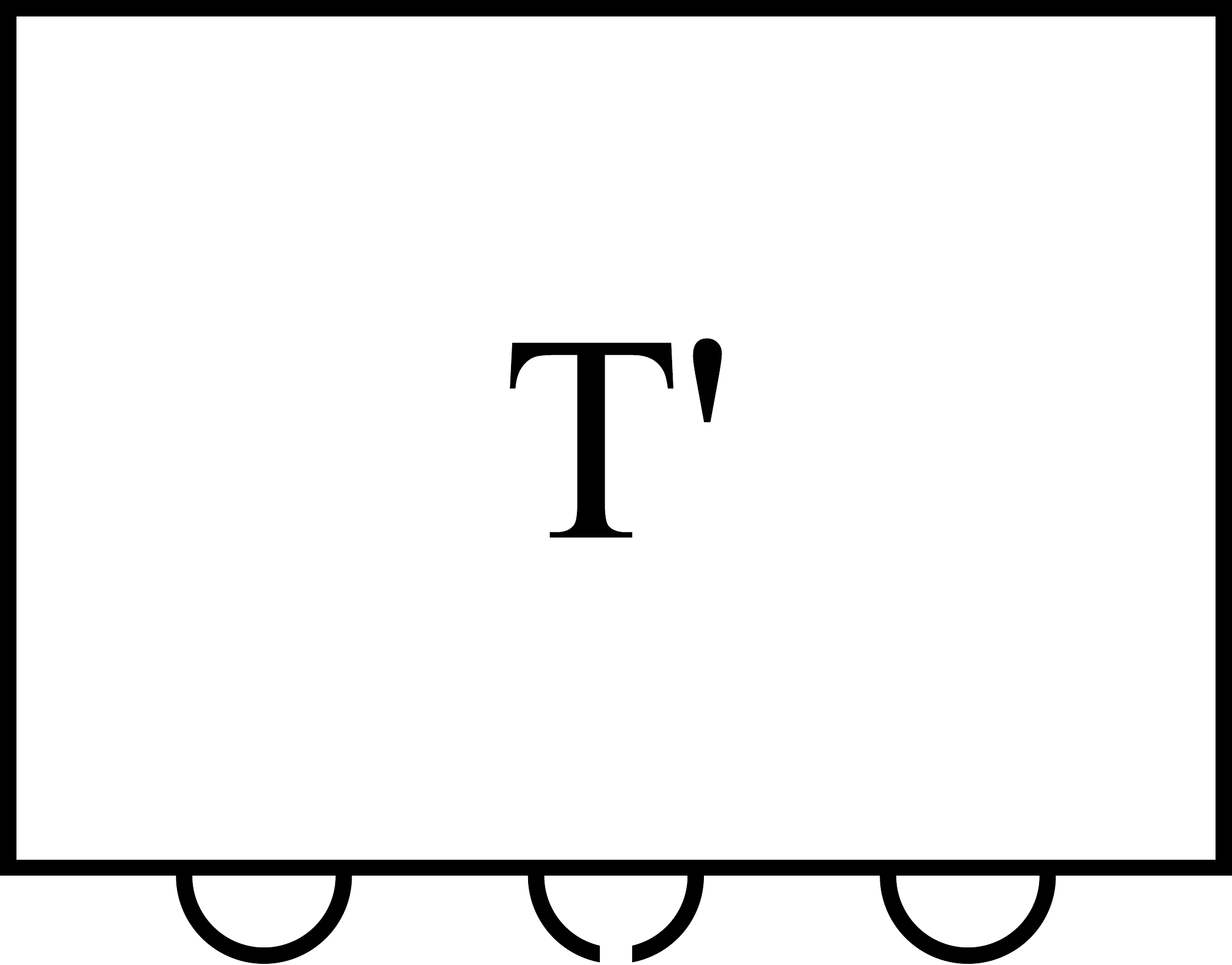}}%
  \end{picture}%
\endgroup%

%% file: diagrams/lemmatwo.pdf_tex
\begingroup%
  \makeatletter%
  \providecommand\color[2][]{%
    \errmessage{(Inkscape) Color is used for the text in Inkscape, but the package 'color.sty' is not loaded}%
    \renewcommand\color[2][]{}%
  }%
  \providecommand\transparent[1]{%
    \errmessage{(Inkscape) Transparency is used (non-zero) for the text in Inkscape, but the package 'transparent.sty' is not loaded}%
    \renewcommand\transparent[1]{}%
  }%
  \providecommand\rotatebox[2]{#2}%
  \ifx\svgwidth\undefined%
    \setlength{\unitlength}{599.27558594bp}%
    \ifx\svgscale\undefined%
      \relax%
    \else%
      \setlength{\unitlength}{\unitlength * \real{\svgscale}}%
    \fi%
  \else%
    \setlength{\unitlength}{\svgwidth}%
  \fi%
  \global\let\svgwidth\undefined%
  \global\let\svgscale\undefined%
  \makeatother%
  \begin{picture}(1,0.71243495)%
    \put(0,0){\includegraphics[width=\unitlength,page=1]{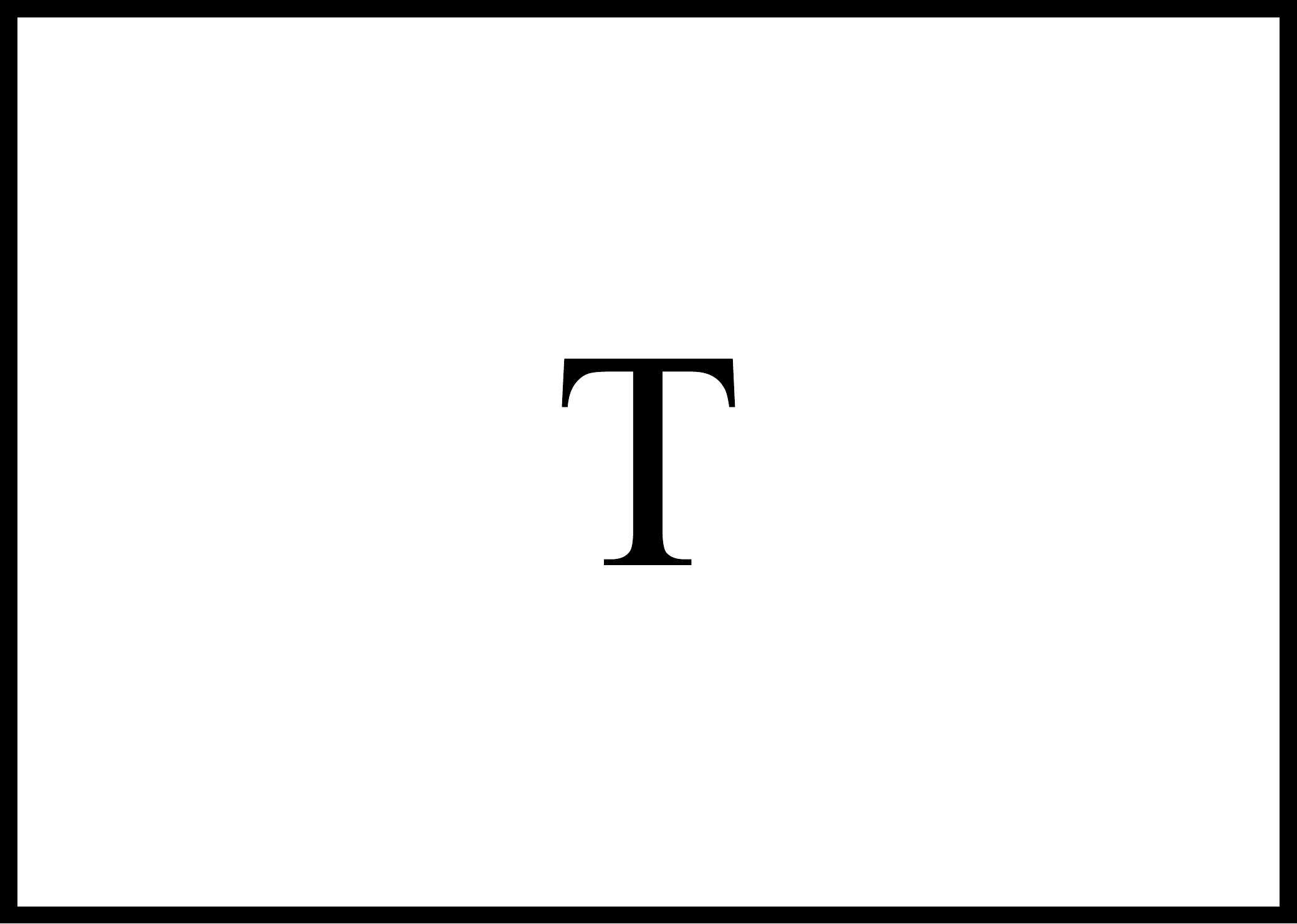}}%
  \end{picture}%
\endgroup%

%% file: diagrams/questionone.pdf_tex
\begingroup%
  \makeatletter%
  \providecommand\color[2][]{%
    \errmessage{(Inkscape) Color is used for the text in Inkscape, but the package 'color.sty' is not loaded}%
    \renewcommand\color[2][]{}%
  }%
  \providecommand\transparent[1]{%
    \errmessage{(Inkscape) Transparency is used (non-zero) for the text in Inkscape, but the package 'transparent.sty' is not loaded}%
    \renewcommand\transparent[1]{}%
  }%
  \providecommand\rotatebox[2]{#2}%
  \ifx\svgwidth\undefined%
    \setlength{\unitlength}{599.27558594bp}%
    \ifx\svgscale\undefined%
      \relax%
    \else%
      \setlength{\unitlength}{\unitlength * \real{\svgscale}}%
    \fi%
  \else%
    \setlength{\unitlength}{\svgwidth}%
  \fi%
  \global\let\svgwidth\undefined%
  \global\let\svgscale\undefined%
  \makeatother%
  \begin{picture}(1,1.30575927)%
    \put(0,0){\includegraphics[width=\unitlength,page=1]{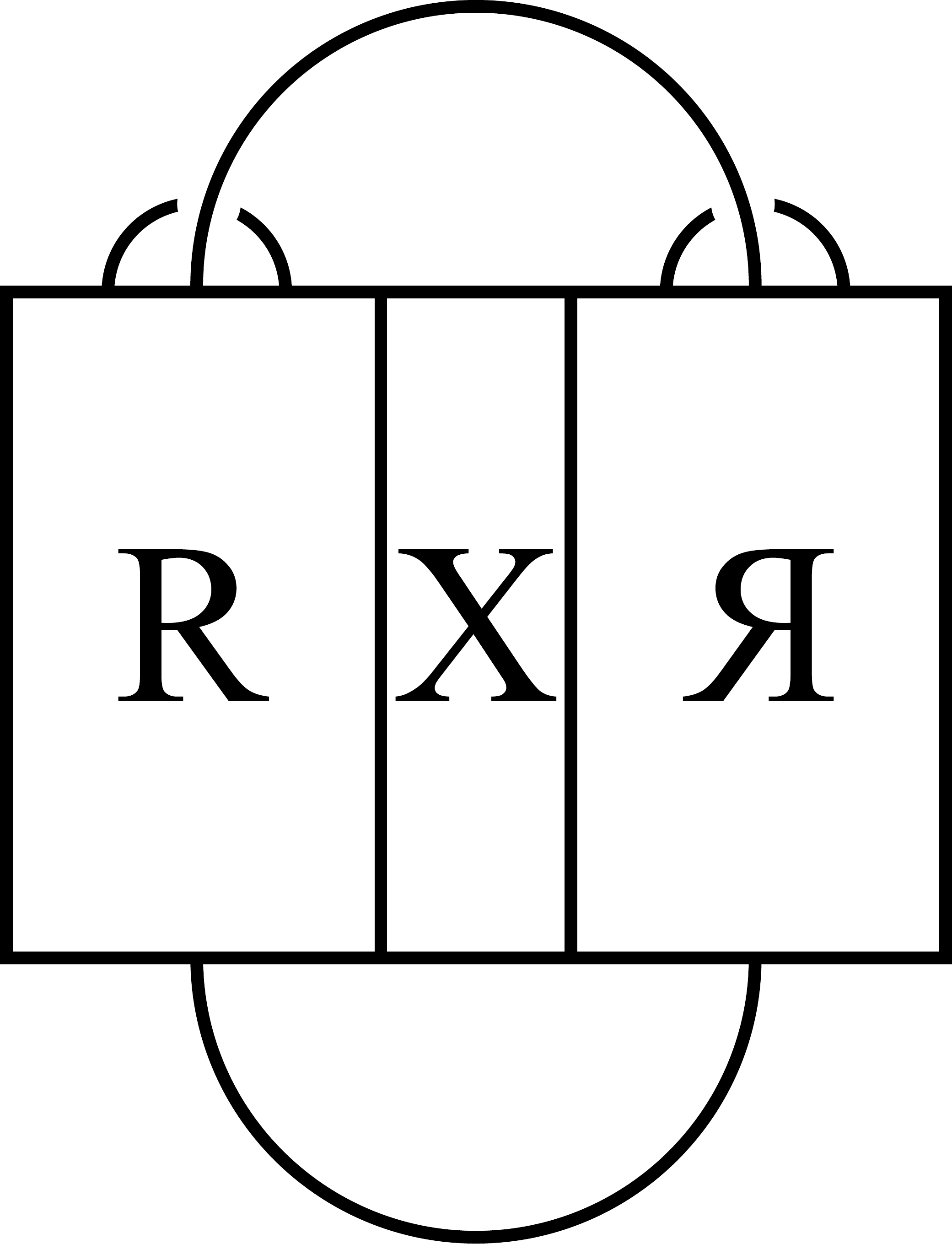}}%
  \end{picture}%
\endgroup%

%% file: diagrams/questiontwo.pdf_tex
\begingroup%
  \makeatletter%
  \providecommand\color[2][]{%
    \errmessage{(Inkscape) Color is used for the text in Inkscape, but the package 'color.sty' is not loaded}%
    \renewcommand\color[2][]{}%
  }%
  \providecommand\transparent[1]{%
    \errmessage{(Inkscape) Transparency is used (non-zero) for the text in Inkscape, but the package 'transparent.sty' is not loaded}%
    \renewcommand\transparent[1]{}%
  }%
  \providecommand\rotatebox[2]{#2}%
  \ifx\svgwidth\undefined%
    \setlength{\unitlength}{599.27558594bp}%
    \ifx\svgscale\undefined%
      \relax%
    \else%
      \setlength{\unitlength}{\unitlength * \real{\svgscale}}%
    \fi%
  \else%
    \setlength{\unitlength}{\svgwidth}%
  \fi%
  \global\let\svgwidth\undefined%
  \global\let\svgscale\undefined%
  \makeatother%
  \begin{picture}(1,1.1057598)%
    \put(0,0){\includegraphics[width=\unitlength,page=1]{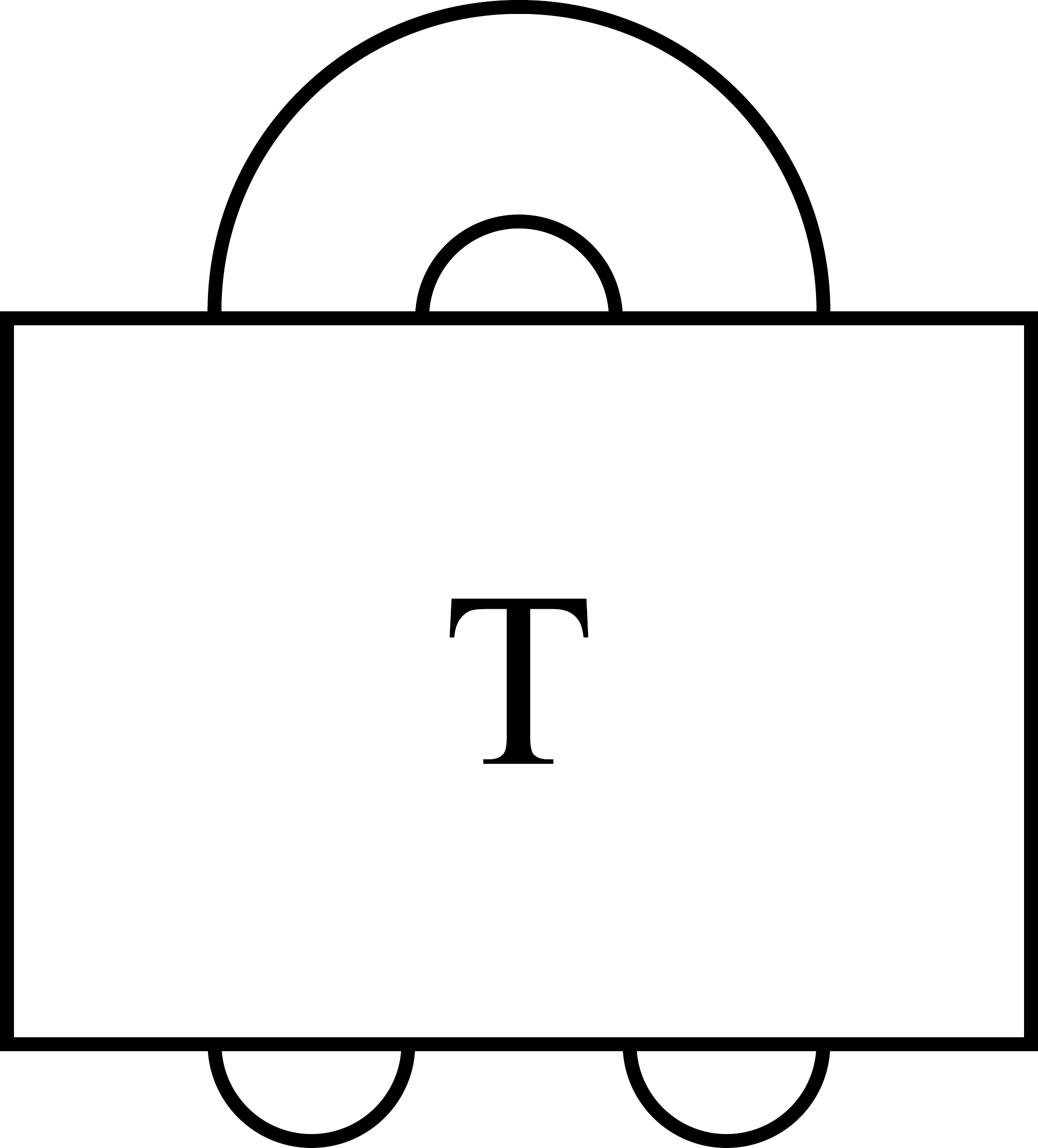}}%
  \end{picture}%
\endgroup%

%% file: diagrams/questionthree.pdf_tex
\begingroup%
  \makeatletter%
  \providecommand\color[2][]{%
    \errmessage{(Inkscape) Color is used for the text in Inkscape, but the package 'color.sty' is not loaded}%
    \renewcommand\color[2][]{}%
  }%
  \providecommand\transparent[1]{%
    \errmessage{(Inkscape) Transparency is used (non-zero) for the text in Inkscape, but the package 'transparent.sty' is not loaded}%
    \renewcommand\transparent[1]{}%
  }%
  \providecommand\rotatebox[2]{#2}%
  \ifx\svgwidth\undefined%
    \setlength{\unitlength}{599.27558594bp}%
    \ifx\svgscale\undefined%
      \relax%
    \else%
      \setlength{\unitlength}{\unitlength * \real{\svgscale}}%
    \fi%
  \else%
    \setlength{\unitlength}{\svgwidth}%
  \fi%
  \global\let\svgwidth\undefined%
  \global\let\svgscale\undefined%
  \makeatother%
  \begin{picture}(1,1.3990839)%
    \put(0,0){\includegraphics[width=\unitlength,page=1]{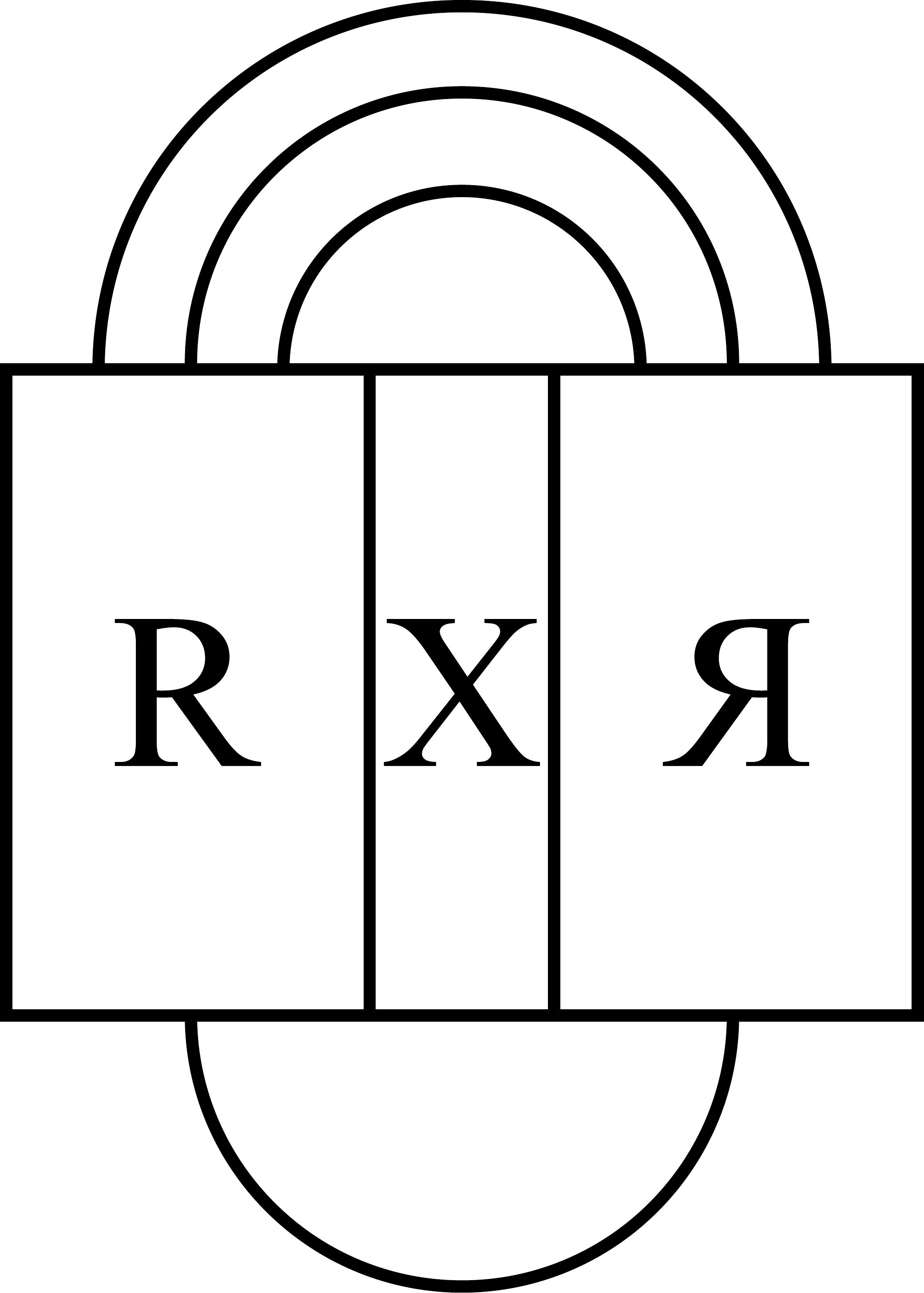}}%
  \end{picture}%
\endgroup%

%% file: diagrams/questionfour.pdf_tex
\begingroup%
  \makeatletter%
  \providecommand\color[2][]{%
    \errmessage{(Inkscape) Color is used for the text in Inkscape, but the package 'color.sty' is not loaded}%
    \renewcommand\color[2][]{}%
  }%
  \providecommand\transparent[1]{%
    \errmessage{(Inkscape) Transparency is used (non-zero) for the text in Inkscape, but the package 'transparent.sty' is not loaded}%
    \renewcommand\transparent[1]{}%
  }%
  \providecommand\rotatebox[2]{#2}%
  \ifx\svgwidth\undefined%
    \setlength{\unitlength}{599.27558594bp}%
    \ifx\svgscale\undefined%
      \relax%
    \else%
      \setlength{\unitlength}{\unitlength * \real{\svgscale}}%
    \fi%
  \else%
    \setlength{\unitlength}{\svgwidth}%
  \fi%
  \global\let\svgwidth\undefined%
  \global\let\svgscale\undefined%
  \makeatother%
  \begin{picture}(1,1.3990839)%
    \put(0,0){\includegraphics[width=\unitlength,page=1]{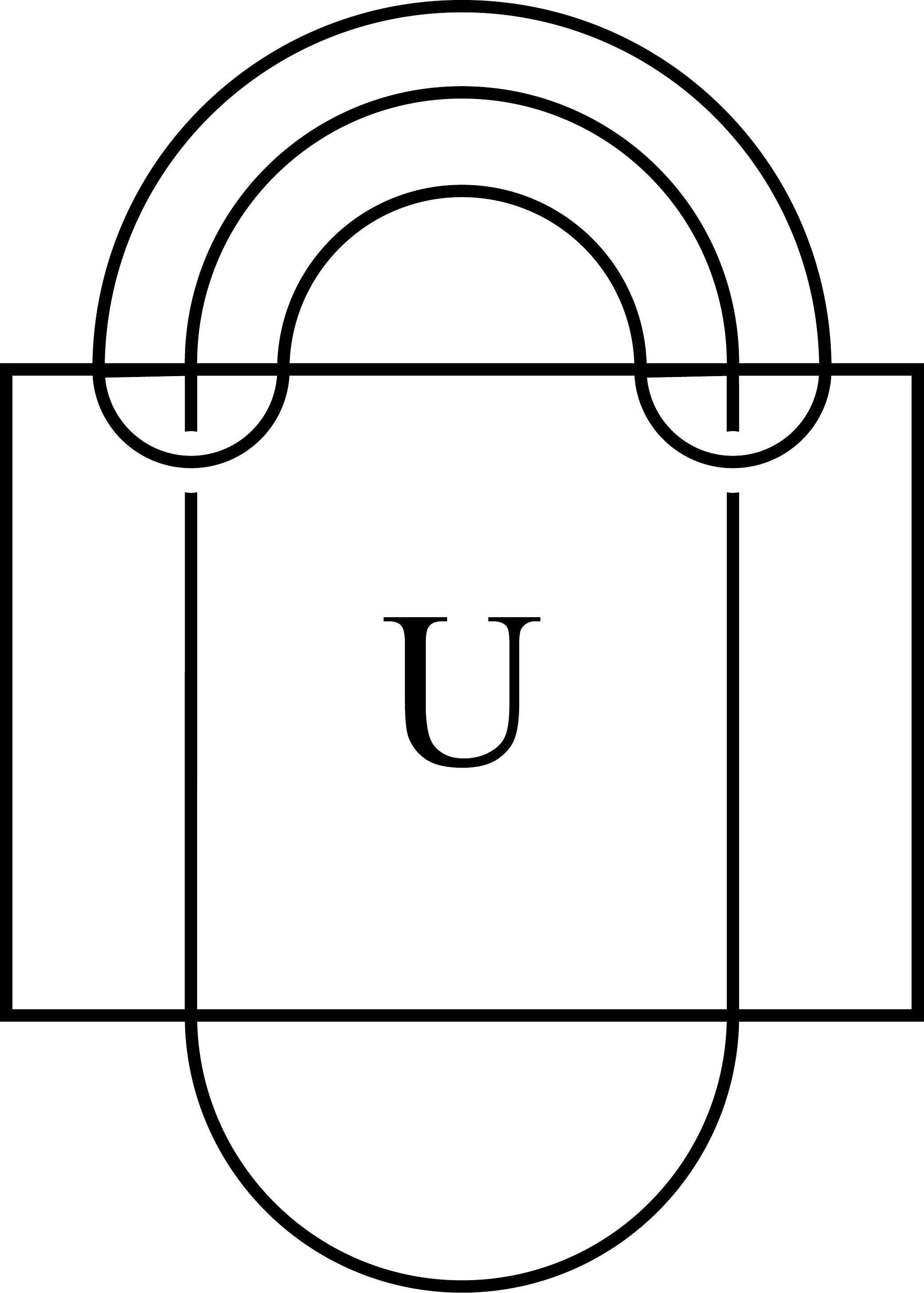}}%
  \end{picture}%
\endgroup%

%% file: diagrams/crossing.pdf_tex
\begingroup%
  \makeatletter%
  \providecommand\color[2][]{%
    \errmessage{(Inkscape) Color is used for the text in Inkscape, but the package 'color.sty' is not loaded}%
    \renewcommand\color[2][]{}%
  }%
  \providecommand\transparent[1]{%
    \errmessage{(Inkscape) Transparency is used (non-zero) for the text in Inkscape, but the package 'transparent.sty' is not loaded}%
    \renewcommand\transparent[1]{}%
  }%
  \providecommand\rotatebox[2]{#2}%
  \ifx\svgwidth\undefined%
    \setlength{\unitlength}{316.54126374bp}%
    \ifx\svgscale\undefined%
      \relax%
    \else%
      \setlength{\unitlength}{\unitlength * \real{\svgscale}}%
    \fi%
  \else%
    \setlength{\unitlength}{\svgwidth}%
  \fi%
  \global\let\svgwidth\undefined%
  \global\let\svgscale\undefined%
  \makeatother%
  \begin{picture}(1,0.98976115)%
    \put(0,0){\includegraphics[width=\unitlength,page=1]{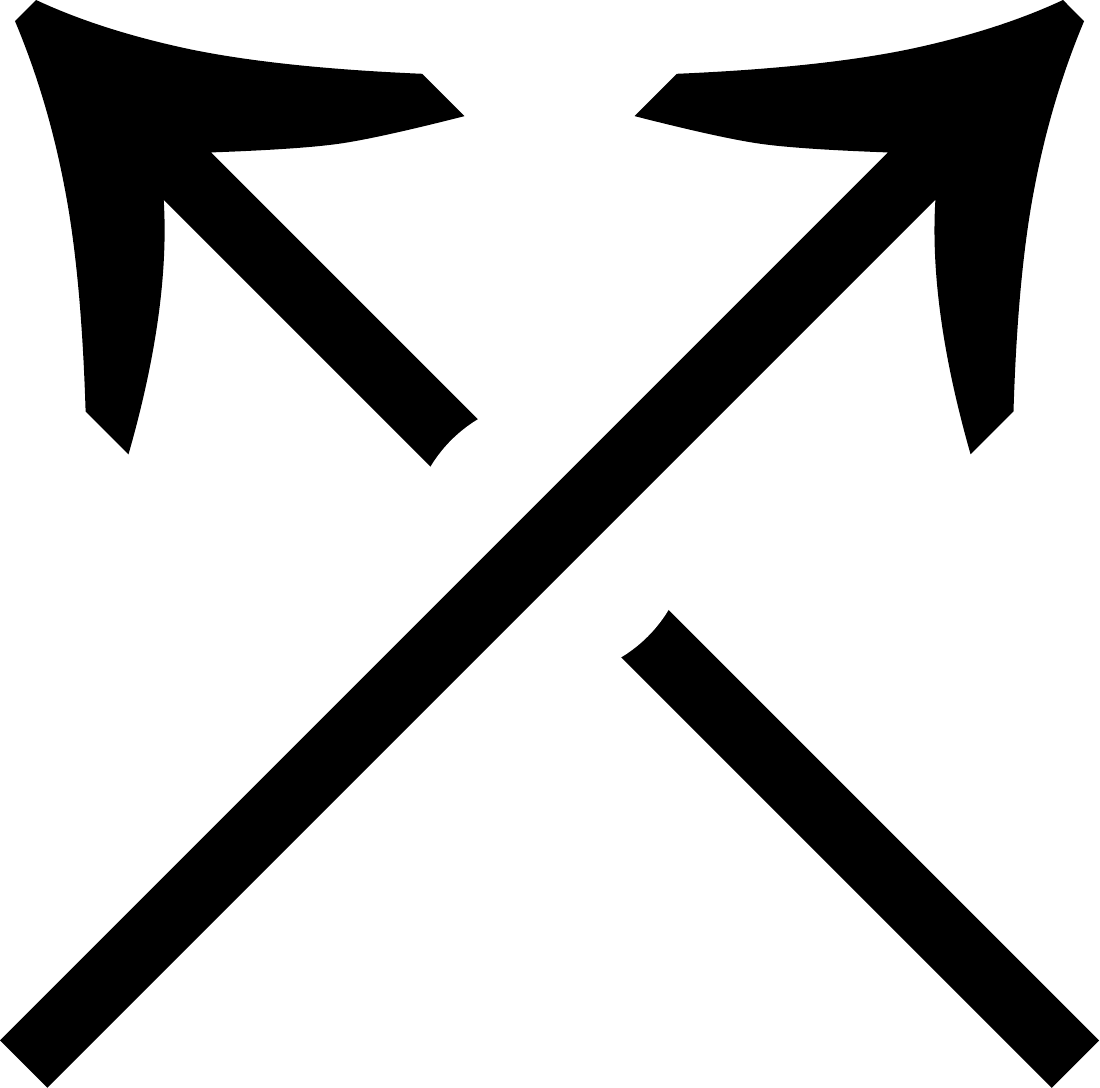}}%
  \end{picture}%
\endgroup%

%% file: diagrams/6_1.pdf_tex
\begingroup%
  \makeatletter%
  \providecommand\color[2][]{%
    \errmessage{(Inkscape) Color is used for the text in Inkscape, but the package 'color.sty' is not loaded}%
    \renewcommand\color[2][]{}%
  }%
  \providecommand\transparent[1]{%
    \errmessage{(Inkscape) Transparency is used (non-zero) for the text in Inkscape, but the package 'transparent.sty' is not loaded}%
    \renewcommand\transparent[1]{}%
  }%
  \providecommand\rotatebox[2]{#2}%
  \ifx\svgwidth\undefined%
    \setlength{\unitlength}{599.275592bp}%
    \ifx\svgscale\undefined%
      \relax%
    \else%
      \setlength{\unitlength}{\unitlength * \real{\svgscale}}%
    \fi%
  \else%
    \setlength{\unitlength}{\svgwidth}%
  \fi%
  \global\let\svgwidth\undefined%
  \global\let\svgscale\undefined%
  \makeatother%
  \begin{picture}(1,1.09985008)%
    \put(0,0){\includegraphics[width=\unitlength,page=1]{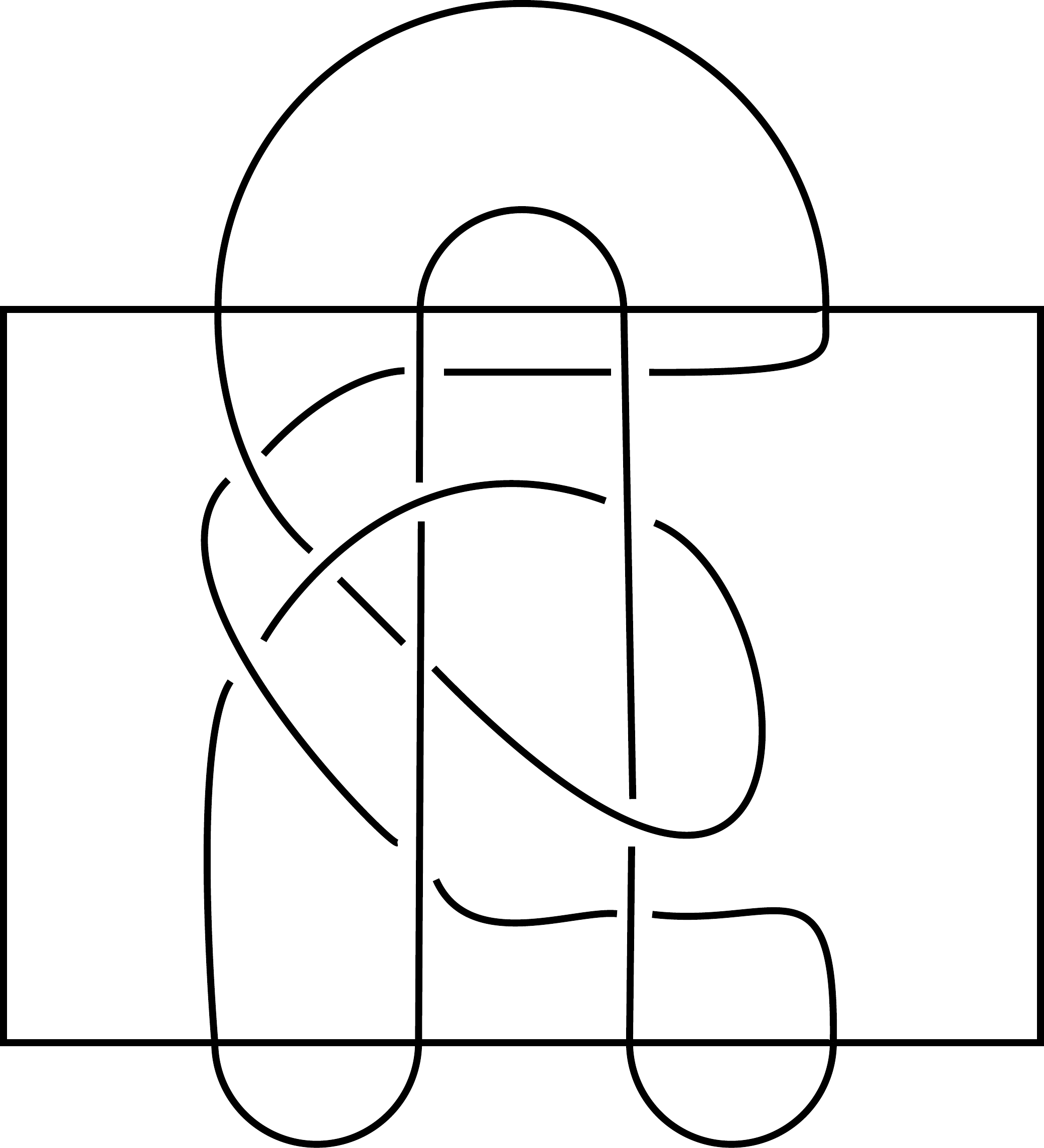}}%
  \end{picture}%
\endgroup%

%% file: diagrams/8_8.pdf_tex
\begingroup%
  \makeatletter%
  \providecommand\color[2][]{%
    \errmessage{(Inkscape) Color is used for the text in Inkscape, but the package 'color.sty' is not loaded}%
    \renewcommand\color[2][]{}%
  }%
  \providecommand\transparent[1]{%
    \errmessage{(Inkscape) Transparency is used (non-zero) for the text in Inkscape, but the package 'transparent.sty' is not loaded}%
    \renewcommand\transparent[1]{}%
  }%
  \providecommand\rotatebox[2]{#2}%
  \ifx\svgwidth\undefined%
    \setlength{\unitlength}{599.275592bp}%
    \ifx\svgscale\undefined%
      \relax%
    \else%
      \setlength{\unitlength}{\unitlength * \real{\svgscale}}%
    \fi%
  \else%
    \setlength{\unitlength}{\svgwidth}%
  \fi%
  \global\let\svgwidth\undefined%
  \global\let\svgscale\undefined%
  \makeatother%
  \begin{picture}(1,1.09985008)%
    \put(0,0){\includegraphics[width=\unitlength,page=1]{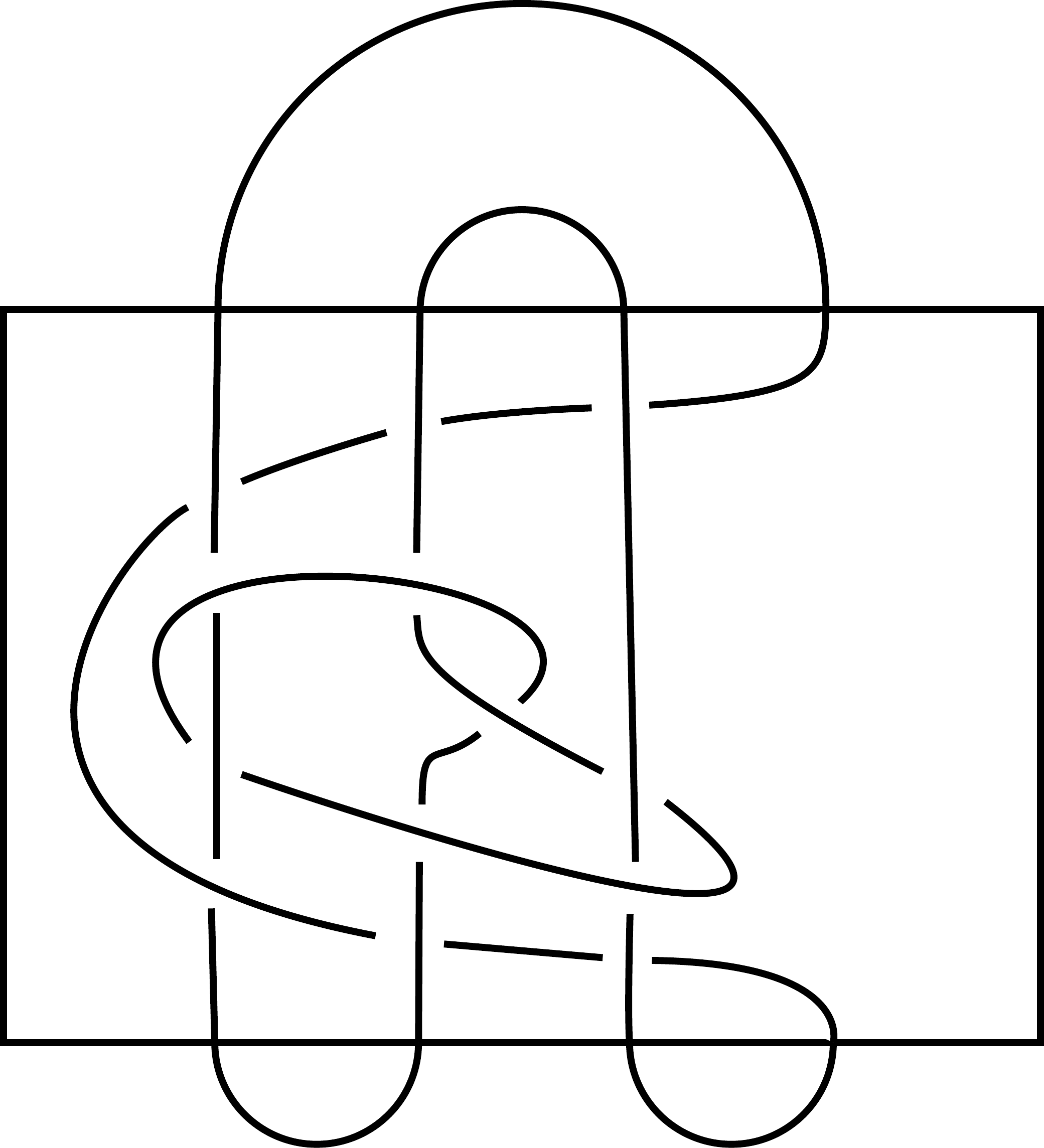}}%
  \end{picture}%
\endgroup%

%% file: diagrams/8_9.pdf_tex
\begingroup%
  \makeatletter%
  \providecommand\color[2][]{%
    \errmessage{(Inkscape) Color is used for the text in Inkscape, but the package 'color.sty' is not loaded}%
    \renewcommand\color[2][]{}%
  }%
  \providecommand\transparent[1]{%
    \errmessage{(Inkscape) Transparency is used (non-zero) for the text in Inkscape, but the package 'transparent.sty' is not loaded}%
    \renewcommand\transparent[1]{}%
  }%
  \providecommand\rotatebox[2]{#2}%
  \ifx\svgwidth\undefined%
    \setlength{\unitlength}{599.275592bp}%
    \ifx\svgscale\undefined%
      \relax%
    \else%
      \setlength{\unitlength}{\unitlength * \real{\svgscale}}%
    \fi%
  \else%
    \setlength{\unitlength}{\svgwidth}%
  \fi%
  \global\let\svgwidth\undefined%
  \global\let\svgscale\undefined%
  \makeatother%
  \begin{picture}(1,1.09985008)%
    \put(0,0){\includegraphics[width=\unitlength,page=1]{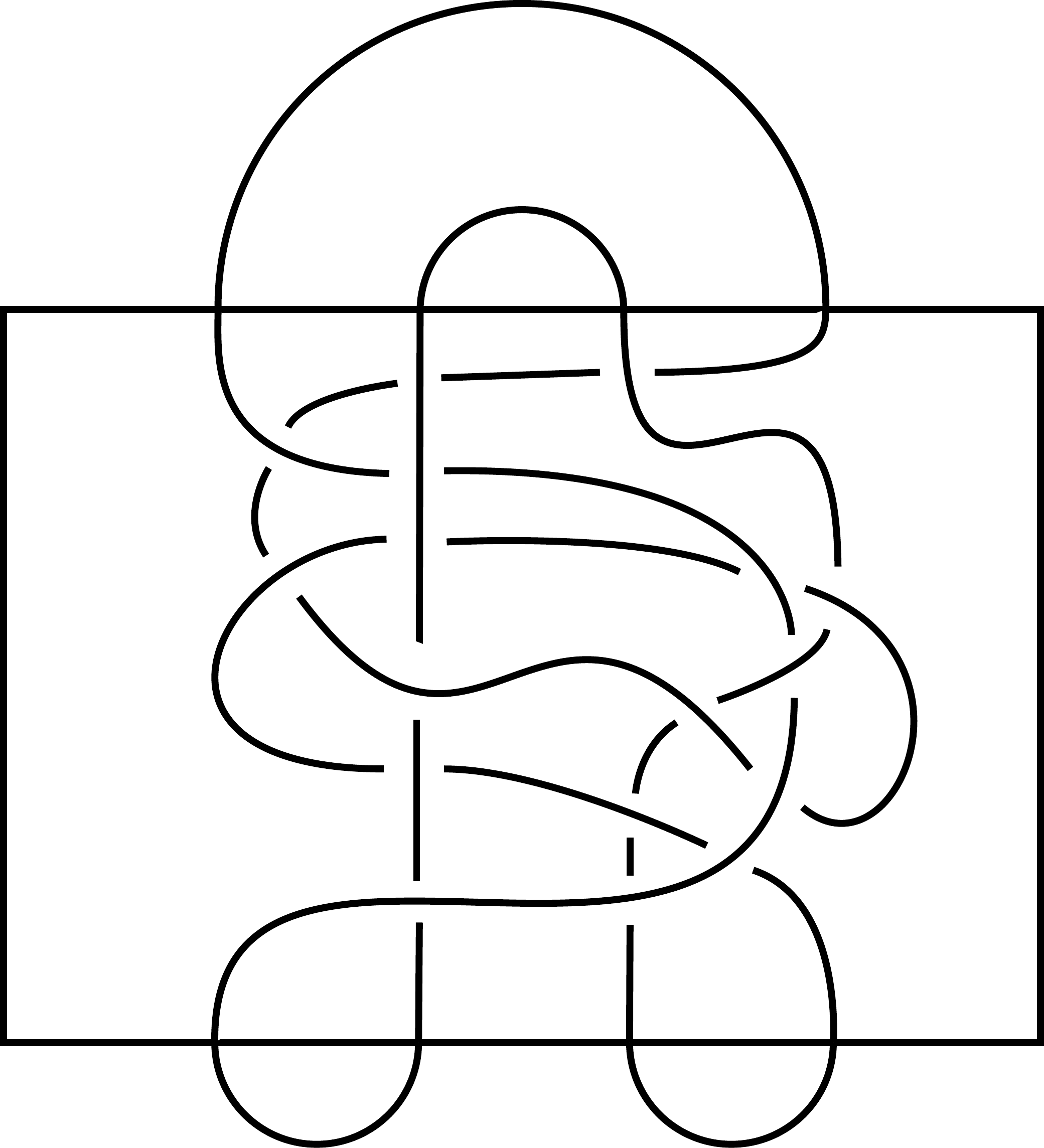}}%
  \end{picture}%
\endgroup%

%% file: diagrams/8_20.pdf_tex
\begingroup%
  \makeatletter%
  \providecommand\color[2][]{%
    \errmessage{(Inkscape) Color is used for the text in Inkscape, but the package 'color.sty' is not loaded}%
    \renewcommand\color[2][]{}%
  }%
  \providecommand\transparent[1]{%
    \errmessage{(Inkscape) Transparency is used (non-zero) for the text in Inkscape, but the package 'transparent.sty' is not loaded}%
    \renewcommand\transparent[1]{}%
  }%
  \providecommand\rotatebox[2]{#2}%
  \ifx\svgwidth\undefined%
    \setlength{\unitlength}{599.275592bp}%
    \ifx\svgscale\undefined%
      \relax%
    \else%
      \setlength{\unitlength}{\unitlength * \real{\svgscale}}%
    \fi%
  \else%
    \setlength{\unitlength}{\svgwidth}%
  \fi%
  \global\let\svgwidth\undefined%
  \global\let\svgscale\undefined%
  \makeatother%
  \begin{picture}(1,1.09985008)%
    \put(0,0){\includegraphics[width=\unitlength,page=1]{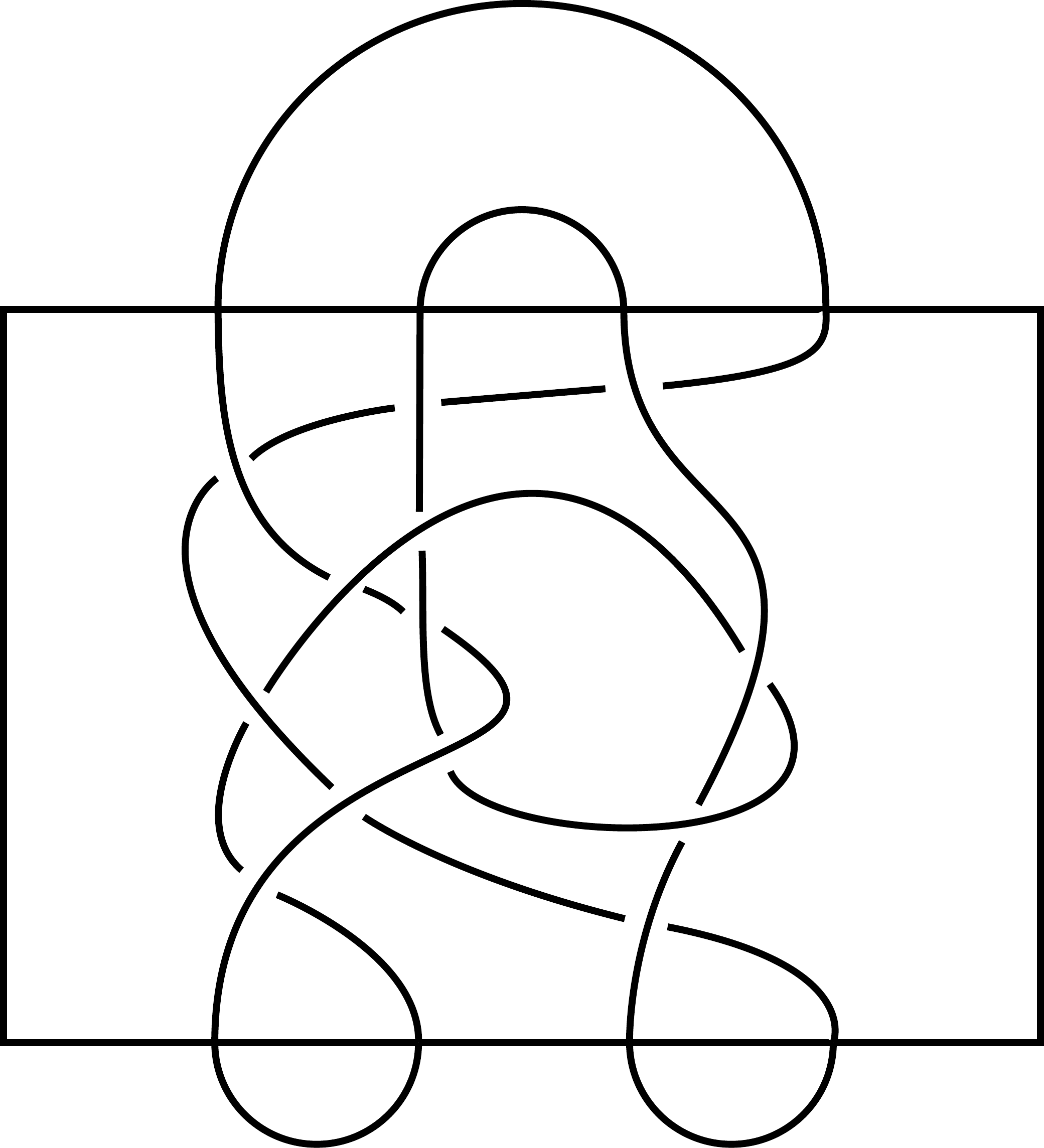}}%
  \end{picture}%
\endgroup%

%% file: diagrams/9_27.pdf_tex
\begingroup%
  \makeatletter%
  \providecommand\color[2][]{%
    \errmessage{(Inkscape) Color is used for the text in Inkscape, but the package 'color.sty' is not loaded}%
    \renewcommand\color[2][]{}%
  }%
  \providecommand\transparent[1]{%
    \errmessage{(Inkscape) Transparency is used (non-zero) for the text in Inkscape, but the package 'transparent.sty' is not loaded}%
    \renewcommand\transparent[1]{}%
  }%
  \providecommand\rotatebox[2]{#2}%
  \ifx\svgwidth\undefined%
    \setlength{\unitlength}{599.275592bp}%
    \ifx\svgscale\undefined%
      \relax%
    \else%
      \setlength{\unitlength}{\unitlength * \real{\svgscale}}%
    \fi%
  \else%
    \setlength{\unitlength}{\svgwidth}%
  \fi%
  \global\let\svgwidth\undefined%
  \global\let\svgscale\undefined%
  \makeatother%
  \begin{picture}(1,1.09985008)%
    \put(0,0){\includegraphics[width=\unitlength,page=1]{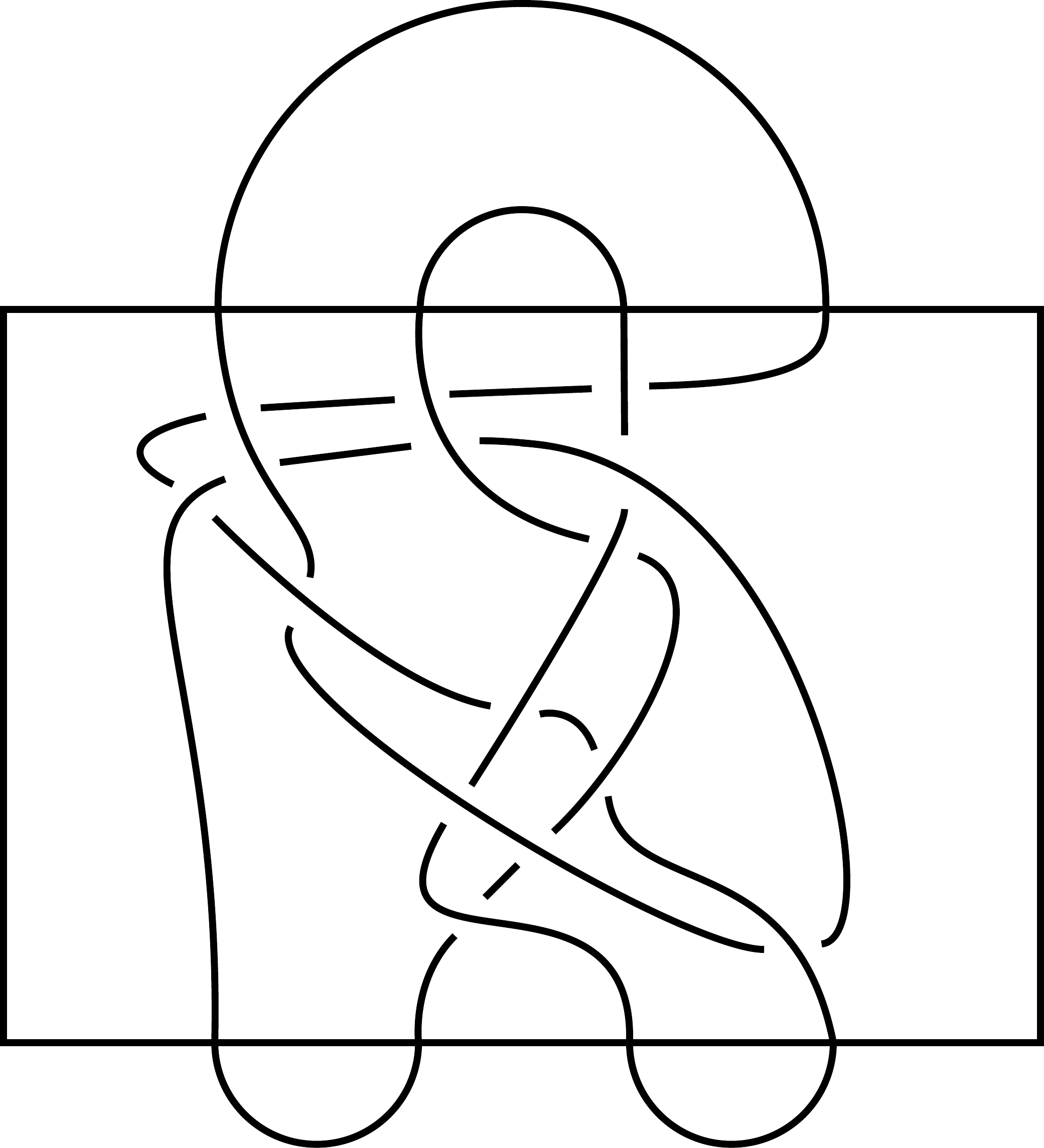}}%
  \end{picture}%
\endgroup%

%% file: diagrams/9_41.pdf_tex
\begingroup%
  \makeatletter%
  \providecommand\color[2][]{%
    \errmessage{(Inkscape) Color is used for the text in Inkscape, but the package 'color.sty' is not loaded}%
    \renewcommand\color[2][]{}%
  }%
  \providecommand\transparent[1]{%
    \errmessage{(Inkscape) Transparency is used (non-zero) for the text in Inkscape, but the package 'transparent.sty' is not loaded}%
    \renewcommand\transparent[1]{}%
  }%
  \providecommand\rotatebox[2]{#2}%
  \ifx\svgwidth\undefined%
    \setlength{\unitlength}{599.275592bp}%
    \ifx\svgscale\undefined%
      \relax%
    \else%
      \setlength{\unitlength}{\unitlength * \real{\svgscale}}%
    \fi%
  \else%
    \setlength{\unitlength}{\svgwidth}%
  \fi%
  \global\let\svgwidth\undefined%
  \global\let\svgscale\undefined%
  \makeatother%
  \begin{picture}(1,1.09985008)%
    \put(0,0){\includegraphics[width=\unitlength,page=1]{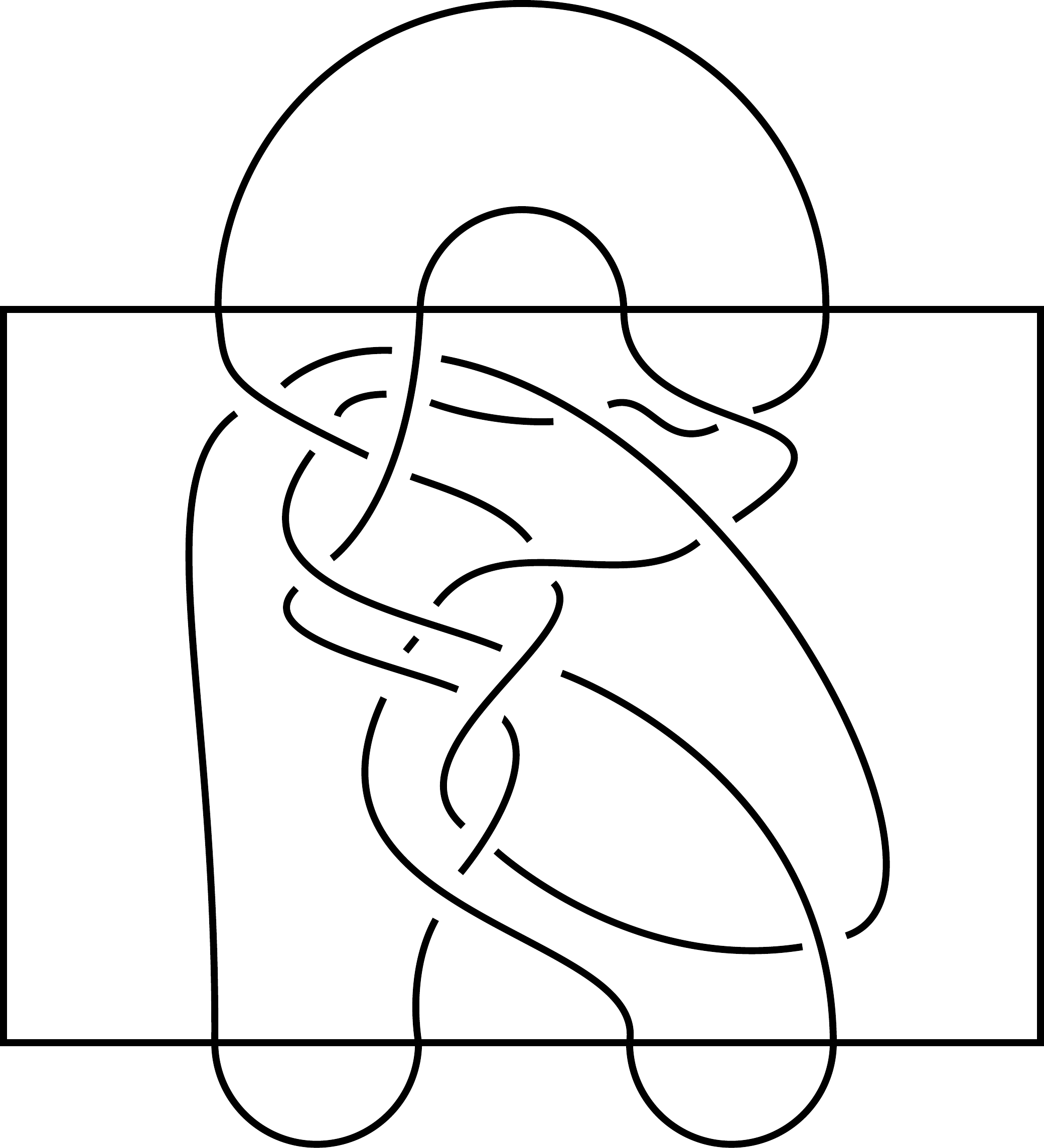}}%
  \end{picture}%
\endgroup%

%% file: diagrams/9_46.pdf_tex
\begingroup%
  \makeatletter%
  \providecommand\color[2][]{%
    \errmessage{(Inkscape) Color is used for the text in Inkscape, but the package 'color.sty' is not loaded}%
    \renewcommand\color[2][]{}%
  }%
  \providecommand\transparent[1]{%
    \errmessage{(Inkscape) Transparency is used (non-zero) for the text in Inkscape, but the package 'transparent.sty' is not loaded}%
    \renewcommand\transparent[1]{}%
  }%
  \providecommand\rotatebox[2]{#2}%
  \ifx\svgwidth\undefined%
    \setlength{\unitlength}{599.275592bp}%
    \ifx\svgscale\undefined%
      \relax%
    \else%
      \setlength{\unitlength}{\unitlength * \real{\svgscale}}%
    \fi%
  \else%
    \setlength{\unitlength}{\svgwidth}%
  \fi%
  \global\let\svgwidth\undefined%
  \global\let\svgscale\undefined%
  \makeatother%
  \begin{picture}(1,1.09985008)%
    \put(0,0){\includegraphics[width=\unitlength,page=1]{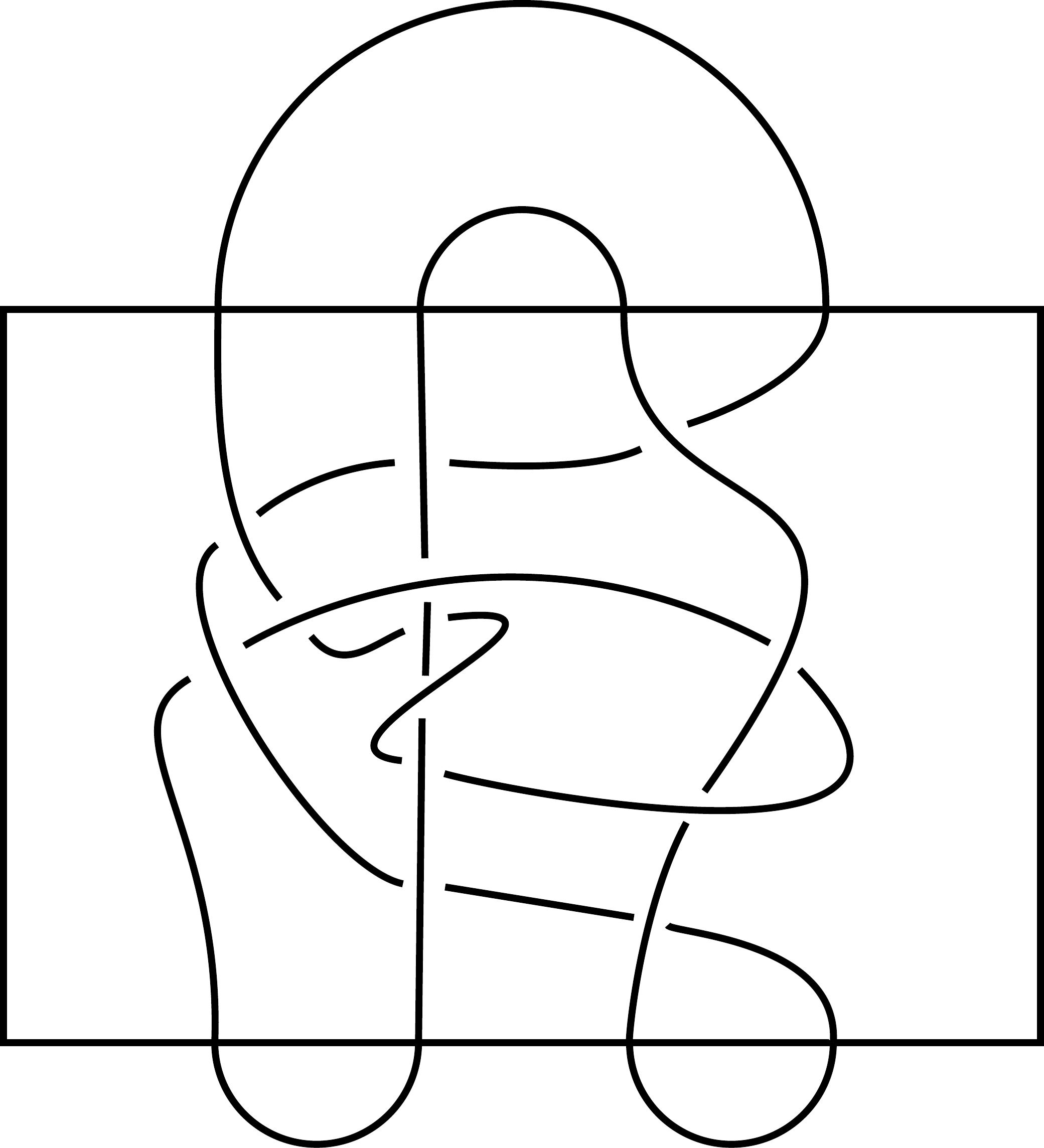}}%
  \end{picture}%
\endgroup%

%% file: diagrams/10_3.pdf_tex
\begingroup%
  \makeatletter%
  \providecommand\color[2][]{%
    \errmessage{(Inkscape) Color is used for the text in Inkscape, but the package 'color.sty' is not loaded}%
    \renewcommand\color[2][]{}%
  }%
  \providecommand\transparent[1]{%
    \errmessage{(Inkscape) Transparency is used (non-zero) for the text in Inkscape, but the package 'transparent.sty' is not loaded}%
    \renewcommand\transparent[1]{}%
  }%
  \providecommand\rotatebox[2]{#2}%
  \ifx\svgwidth\undefined%
    \setlength{\unitlength}{599.275592bp}%
    \ifx\svgscale\undefined%
      \relax%
    \else%
      \setlength{\unitlength}{\unitlength * \real{\svgscale}}%
    \fi%
  \else%
    \setlength{\unitlength}{\svgwidth}%
  \fi%
  \global\let\svgwidth\undefined%
  \global\let\svgscale\undefined%
  \makeatother%
  \begin{picture}(1,1.09985008)%
    \put(0,0){\includegraphics[width=\unitlength,page=1]{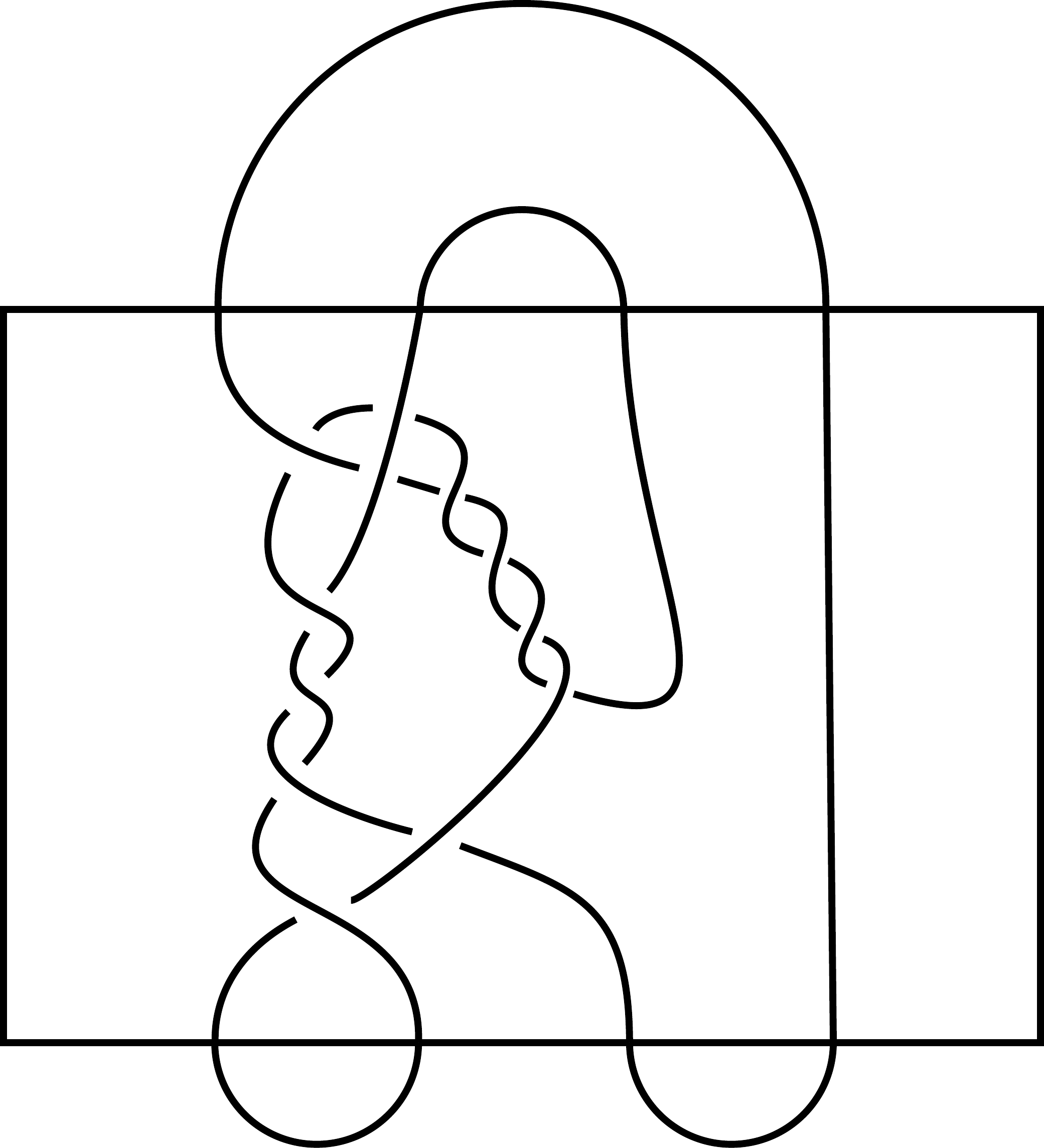}}%
  \end{picture}%
\endgroup%

%% file: diagrams/10_22.pdf_tex
\begingroup%
  \makeatletter%
  \providecommand\color[2][]{%
    \errmessage{(Inkscape) Color is used for the text in Inkscape, but the package 'color.sty' is not loaded}%
    \renewcommand\color[2][]{}%
  }%
  \providecommand\transparent[1]{%
    \errmessage{(Inkscape) Transparency is used (non-zero) for the text in Inkscape, but the package 'transparent.sty' is not loaded}%
    \renewcommand\transparent[1]{}%
  }%
  \providecommand\rotatebox[2]{#2}%
  \ifx\svgwidth\undefined%
    \setlength{\unitlength}{599.275592bp}%
    \ifx\svgscale\undefined%
      \relax%
    \else%
      \setlength{\unitlength}{\unitlength * \real{\svgscale}}%
    \fi%
  \else%
    \setlength{\unitlength}{\svgwidth}%
  \fi%
  \global\let\svgwidth\undefined%
  \global\let\svgscale\undefined%
  \makeatother%
  \begin{picture}(1,1.09985008)%
    \put(0,0){\includegraphics[width=\unitlength,page=1]{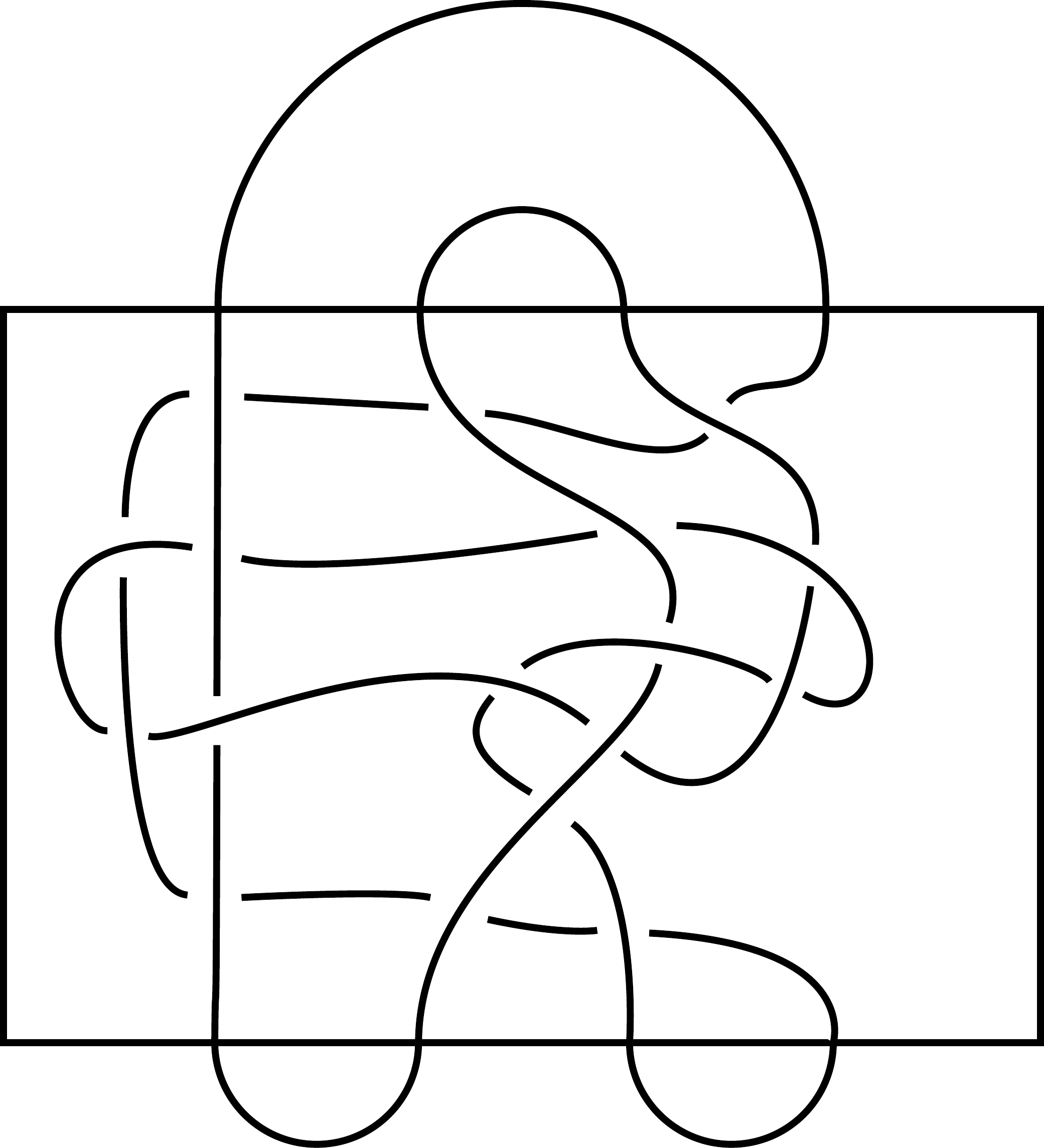}}%
  \end{picture}%
\endgroup%

%% file: diagrams/10_35.pdf_tex
\begingroup%
  \makeatletter%
  \providecommand\color[2][]{%
    \errmessage{(Inkscape) Color is used for the text in Inkscape, but the package 'color.sty' is not loaded}%
    \renewcommand\color[2][]{}%
  }%
  \providecommand\transparent[1]{%
    \errmessage{(Inkscape) Transparency is used (non-zero) for the text in Inkscape, but the package 'transparent.sty' is not loaded}%
    \renewcommand\transparent[1]{}%
  }%
  \providecommand\rotatebox[2]{#2}%
  \ifx\svgwidth\undefined%
    \setlength{\unitlength}{599.275592bp}%
    \ifx\svgscale\undefined%
      \relax%
    \else%
      \setlength{\unitlength}{\unitlength * \real{\svgscale}}%
    \fi%
  \else%
    \setlength{\unitlength}{\svgwidth}%
  \fi%
  \global\let\svgwidth\undefined%
  \global\let\svgscale\undefined%
  \makeatother%
  \begin{picture}(1,1.09985008)%
    \put(0,0){\includegraphics[width=\unitlength,page=1]{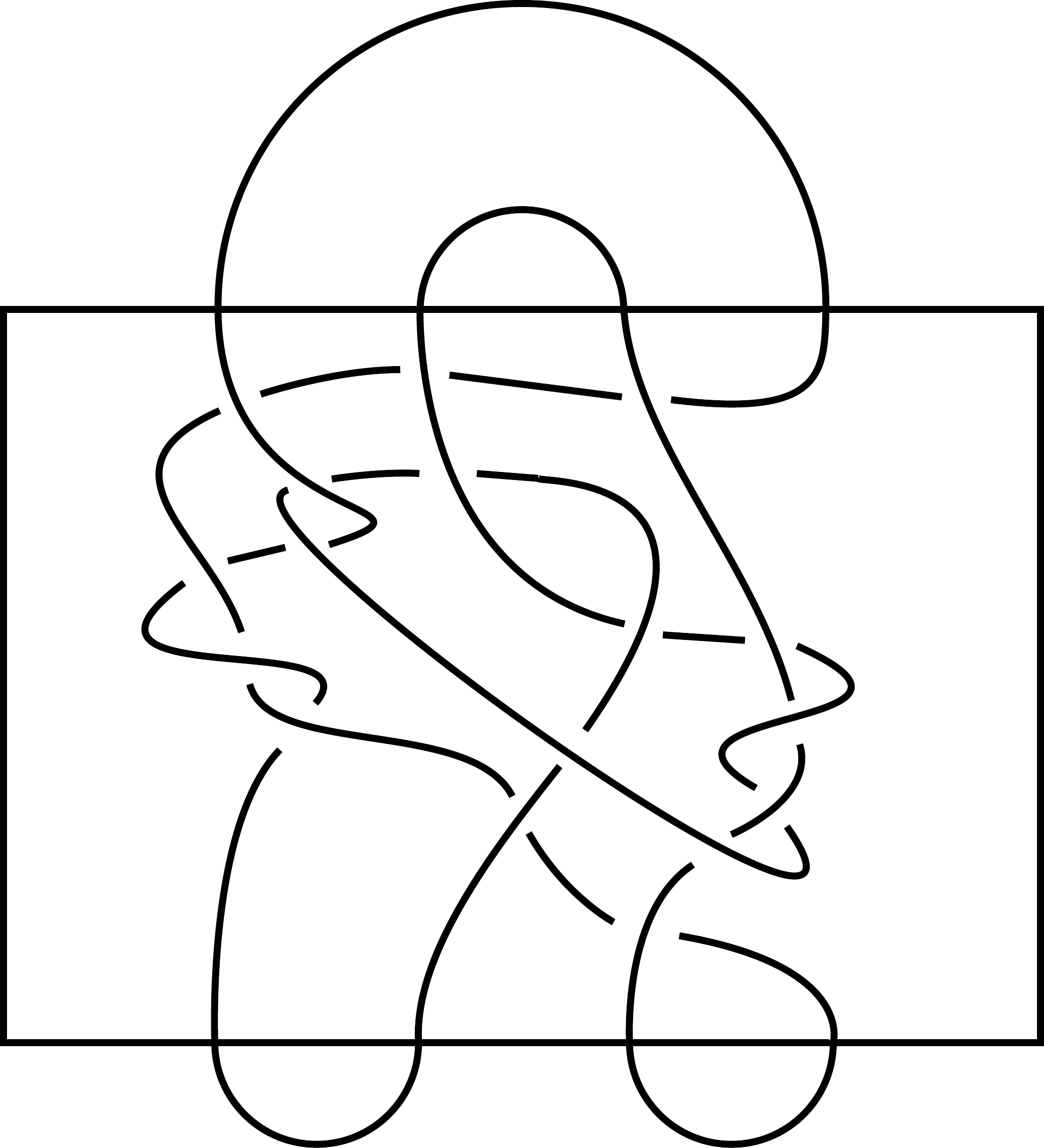}}%
  \end{picture}%
\endgroup%

%% file: diagrams/10_42.pdf_tex
\begingroup%
  \makeatletter%
  \providecommand\color[2][]{%
    \errmessage{(Inkscape) Color is used for the text in Inkscape, but the package 'color.sty' is not loaded}%
    \renewcommand\color[2][]{}%
  }%
  \providecommand\transparent[1]{%
    \errmessage{(Inkscape) Transparency is used (non-zero) for the text in Inkscape, but the package 'transparent.sty' is not loaded}%
    \renewcommand\transparent[1]{}%
  }%
  \providecommand\rotatebox[2]{#2}%
  \ifx\svgwidth\undefined%
    \setlength{\unitlength}{599.275592bp}%
    \ifx\svgscale\undefined%
      \relax%
    \else%
      \setlength{\unitlength}{\unitlength * \real{\svgscale}}%
    \fi%
  \else%
    \setlength{\unitlength}{\svgwidth}%
  \fi%
  \global\let\svgwidth\undefined%
  \global\let\svgscale\undefined%
  \makeatother%
  \begin{picture}(1,1.09985008)%
    \put(0,0){\includegraphics[width=\unitlength,page=1]{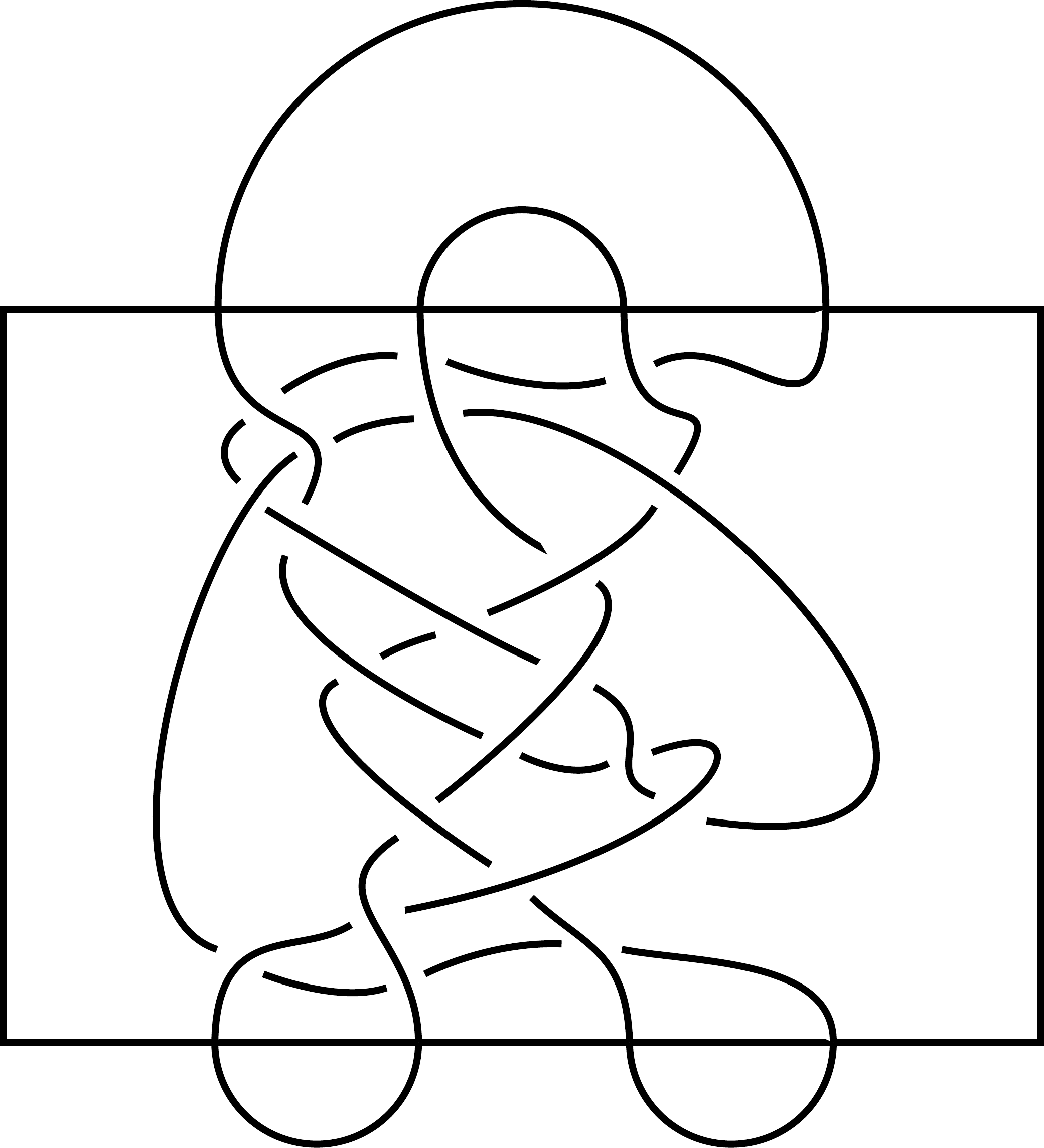}}%
  \end{picture}%
\endgroup%

%% file: diagrams/10_48.pdf_tex
\begingroup%
  \makeatletter%
  \providecommand\color[2][]{%
    \errmessage{(Inkscape) Color is used for the text in Inkscape, but the package 'color.sty' is not loaded}%
    \renewcommand\color[2][]{}%
  }%
  \providecommand\transparent[1]{%
    \errmessage{(Inkscape) Transparency is used (non-zero) for the text in Inkscape, but the package 'transparent.sty' is not loaded}%
    \renewcommand\transparent[1]{}%
  }%
  \providecommand\rotatebox[2]{#2}%
  \ifx\svgwidth\undefined%
    \setlength{\unitlength}{599.275592bp}%
    \ifx\svgscale\undefined%
      \relax%
    \else%
      \setlength{\unitlength}{\unitlength * \real{\svgscale}}%
    \fi%
  \else%
    \setlength{\unitlength}{\svgwidth}%
  \fi%
  \global\let\svgwidth\undefined%
  \global\let\svgscale\undefined%
  \makeatother%
  \begin{picture}(1,1.09985008)%
    \put(0,0){\includegraphics[width=\unitlength,page=1]{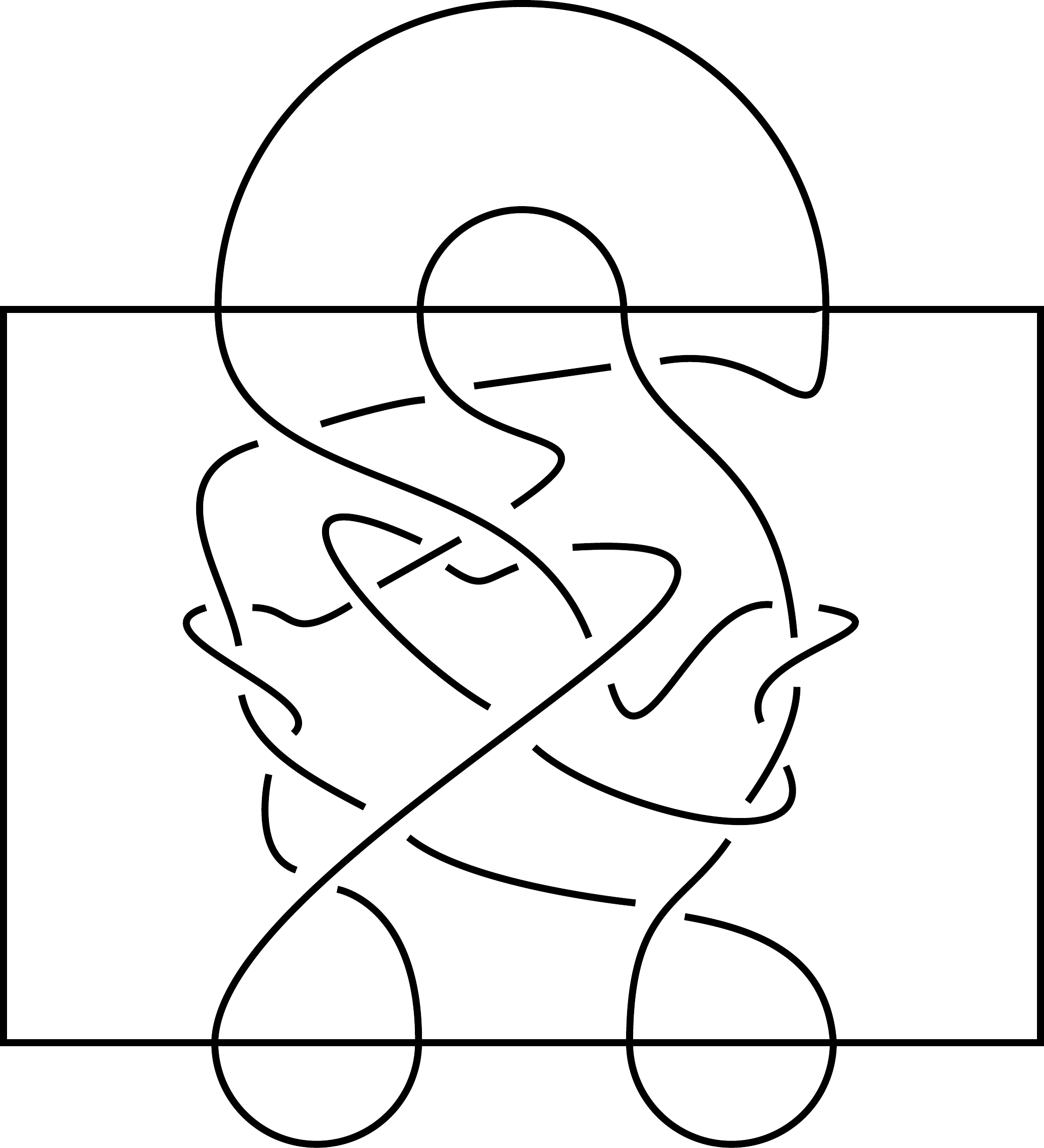}}%
  \end{picture}%
\endgroup%

%% file: diagrams/10_75.pdf_tex
\begingroup%
  \makeatletter%
  \providecommand\color[2][]{%
    \errmessage{(Inkscape) Color is used for the text in Inkscape, but the package 'color.sty' is not loaded}%
    \renewcommand\color[2][]{}%
  }%
  \providecommand\transparent[1]{%
    \errmessage{(Inkscape) Transparency is used (non-zero) for the text in Inkscape, but the package 'transparent.sty' is not loaded}%
    \renewcommand\transparent[1]{}%
  }%
  \providecommand\rotatebox[2]{#2}%
  \ifx\svgwidth\undefined%
    \setlength{\unitlength}{599.27558594bp}%
    \ifx\svgscale\undefined%
      \relax%
    \else%
      \setlength{\unitlength}{\unitlength * \real{\svgscale}}%
    \fi%
  \else%
    \setlength{\unitlength}{\svgwidth}%
  \fi%
  \global\let\svgwidth\undefined%
  \global\let\svgscale\undefined%
  \makeatother%
  \begin{picture}(1,1.09985009)%
    \put(0,0){\includegraphics[width=\unitlength,page=1]{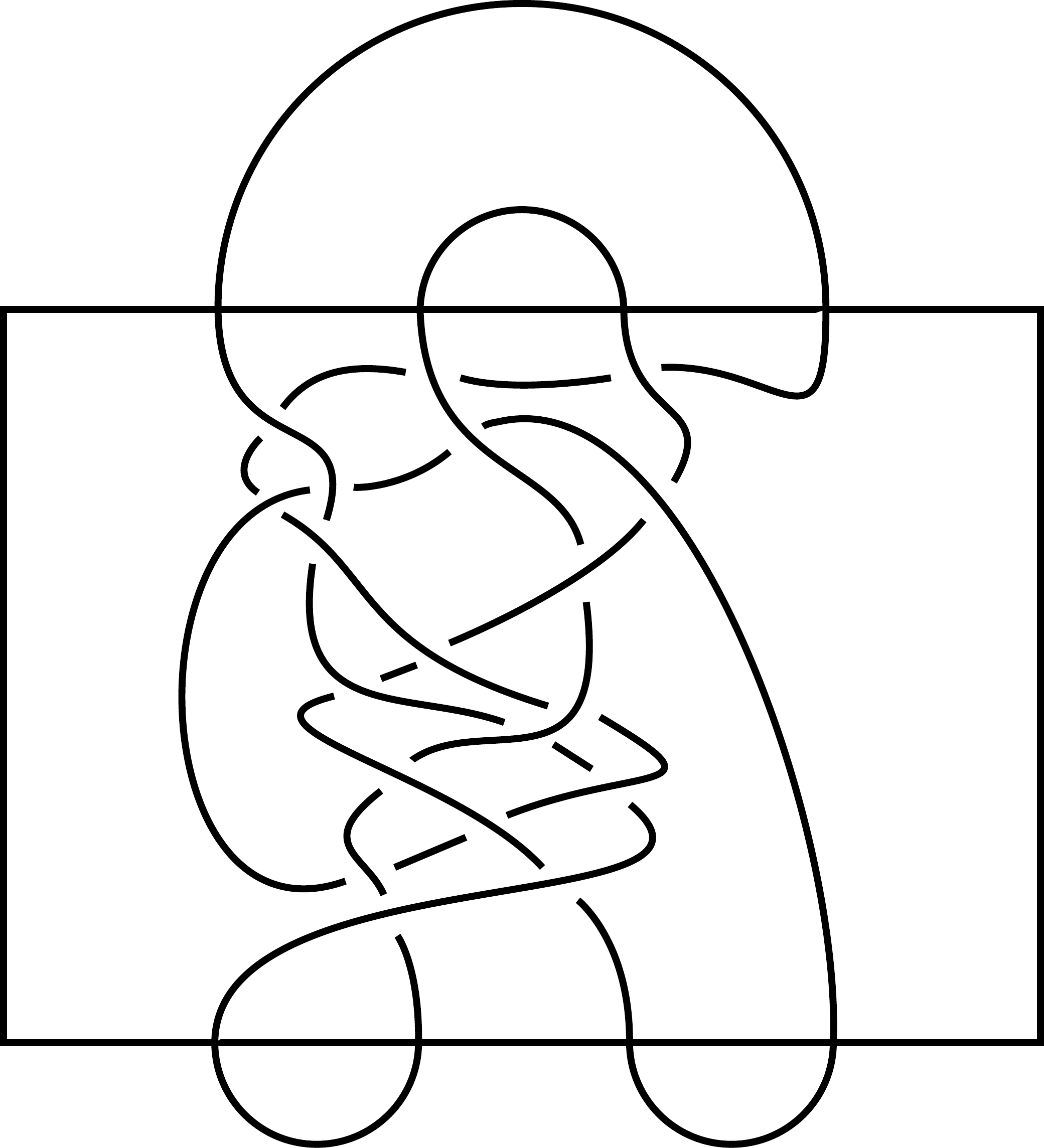}}%
  \end{picture}%
\endgroup%

%% file: diagrams/10_87.pdf_tex
\begingroup%
  \makeatletter%
  \providecommand\color[2][]{%
    \errmessage{(Inkscape) Color is used for the text in Inkscape, but the package 'color.sty' is not loaded}%
    \renewcommand\color[2][]{}%
  }%
  \providecommand\transparent[1]{%
    \errmessage{(Inkscape) Transparency is used (non-zero) for the text in Inkscape, but the package 'transparent.sty' is not loaded}%
    \renewcommand\transparent[1]{}%
  }%
  \providecommand\rotatebox[2]{#2}%
  \ifx\svgwidth\undefined%
    \setlength{\unitlength}{599.275592bp}%
    \ifx\svgscale\undefined%
      \relax%
    \else%
      \setlength{\unitlength}{\unitlength * \real{\svgscale}}%
    \fi%
  \else%
    \setlength{\unitlength}{\svgwidth}%
  \fi%
  \global\let\svgwidth\undefined%
  \global\let\svgscale\undefined%
  \makeatother%
  \begin{picture}(1,1.09985008)%
    \put(0,0){\includegraphics[width=\unitlength,page=1]{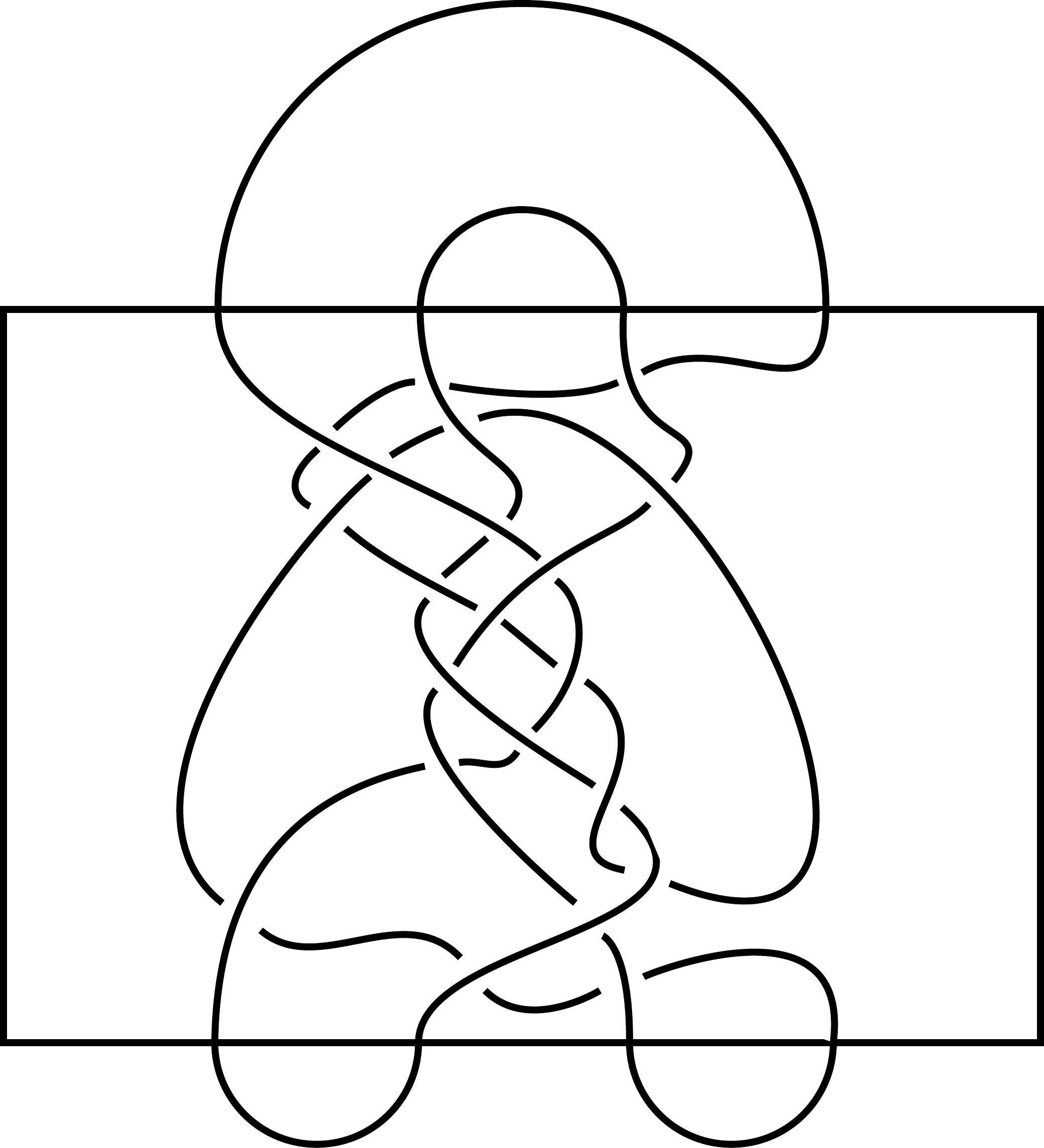}}%
  \end{picture}%
\endgroup%

%% file: diagrams/10_99.pdf_tex
\begingroup%
  \makeatletter%
  \providecommand\color[2][]{%
    \errmessage{(Inkscape) Color is used for the text in Inkscape, but the package 'color.sty' is not loaded}%
    \renewcommand\color[2][]{}%
  }%
  \providecommand\transparent[1]{%
    \errmessage{(Inkscape) Transparency is used (non-zero) for the text in Inkscape, but the package 'transparent.sty' is not loaded}%
    \renewcommand\transparent[1]{}%
  }%
  \providecommand\rotatebox[2]{#2}%
  \ifx\svgwidth\undefined%
    \setlength{\unitlength}{599.275592bp}%
    \ifx\svgscale\undefined%
      \relax%
    \else%
      \setlength{\unitlength}{\unitlength * \real{\svgscale}}%
    \fi%
  \else%
    \setlength{\unitlength}{\svgwidth}%
  \fi%
  \global\let\svgwidth\undefined%
  \global\let\svgscale\undefined%
  \makeatother%
  \begin{picture}(1,1.09985008)%
    \put(0,0){\includegraphics[width=\unitlength,page=1]{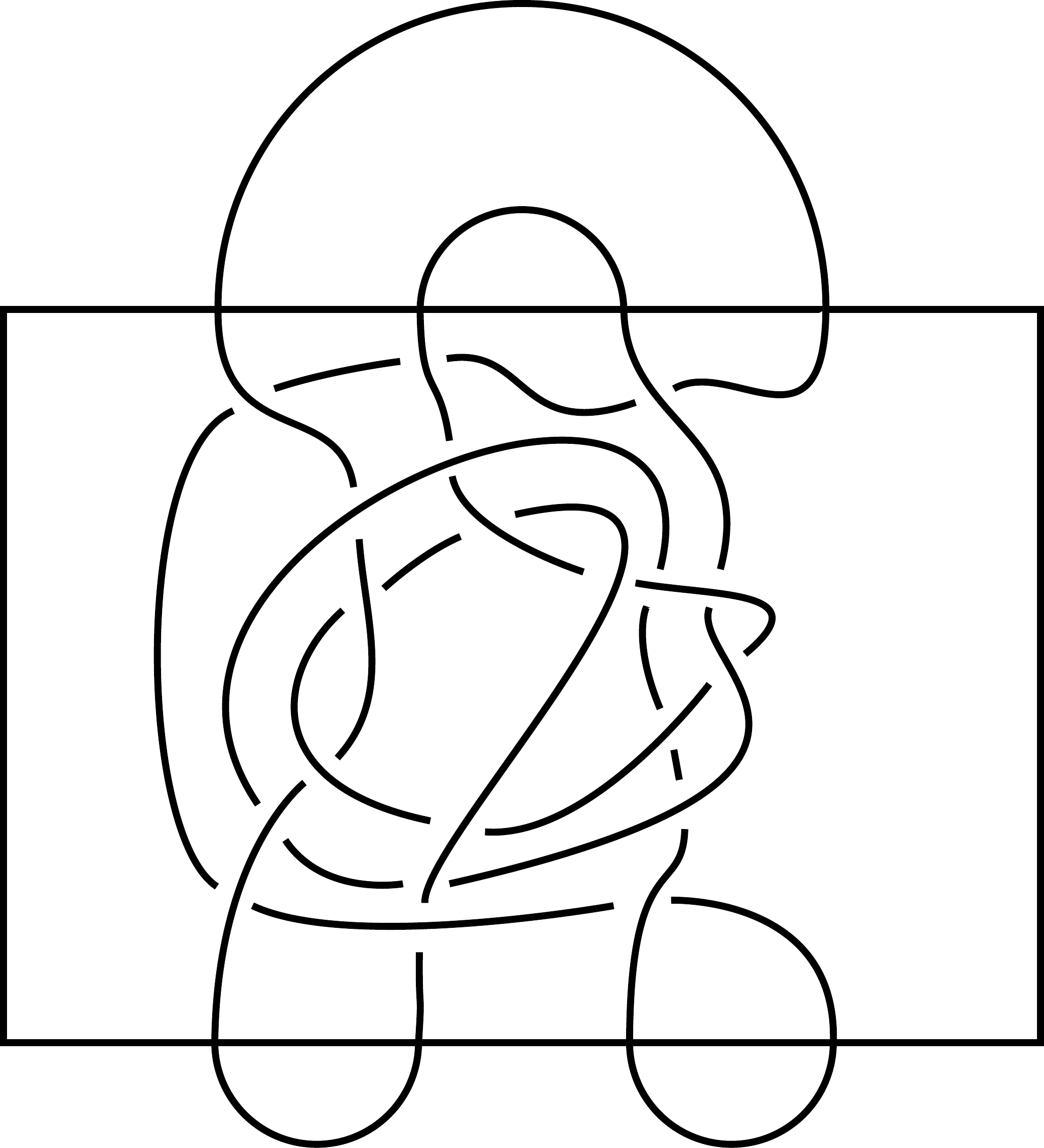}}%
  \end{picture}%
\endgroup%

%% file: diagrams/10_123.pdf_tex
\begingroup%
  \makeatletter%
  \providecommand\color[2][]{%
    \errmessage{(Inkscape) Color is used for the text in Inkscape, but the package 'color.sty' is not loaded}%
    \renewcommand\color[2][]{}%
  }%
  \providecommand\transparent[1]{%
    \errmessage{(Inkscape) Transparency is used (non-zero) for the text in Inkscape, but the package 'transparent.sty' is not loaded}%
    \renewcommand\transparent[1]{}%
  }%
  \providecommand\rotatebox[2]{#2}%
  \ifx\svgwidth\undefined%
    \setlength{\unitlength}{599.275592bp}%
    \ifx\svgscale\undefined%
      \relax%
    \else%
      \setlength{\unitlength}{\unitlength * \real{\svgscale}}%
    \fi%
  \else%
    \setlength{\unitlength}{\svgwidth}%
  \fi%
  \global\let\svgwidth\undefined%
  \global\let\svgscale\undefined%
  \makeatother%
  \begin{picture}(1,1.09985008)%
    \put(0,0){\includegraphics[width=\unitlength,page=1]{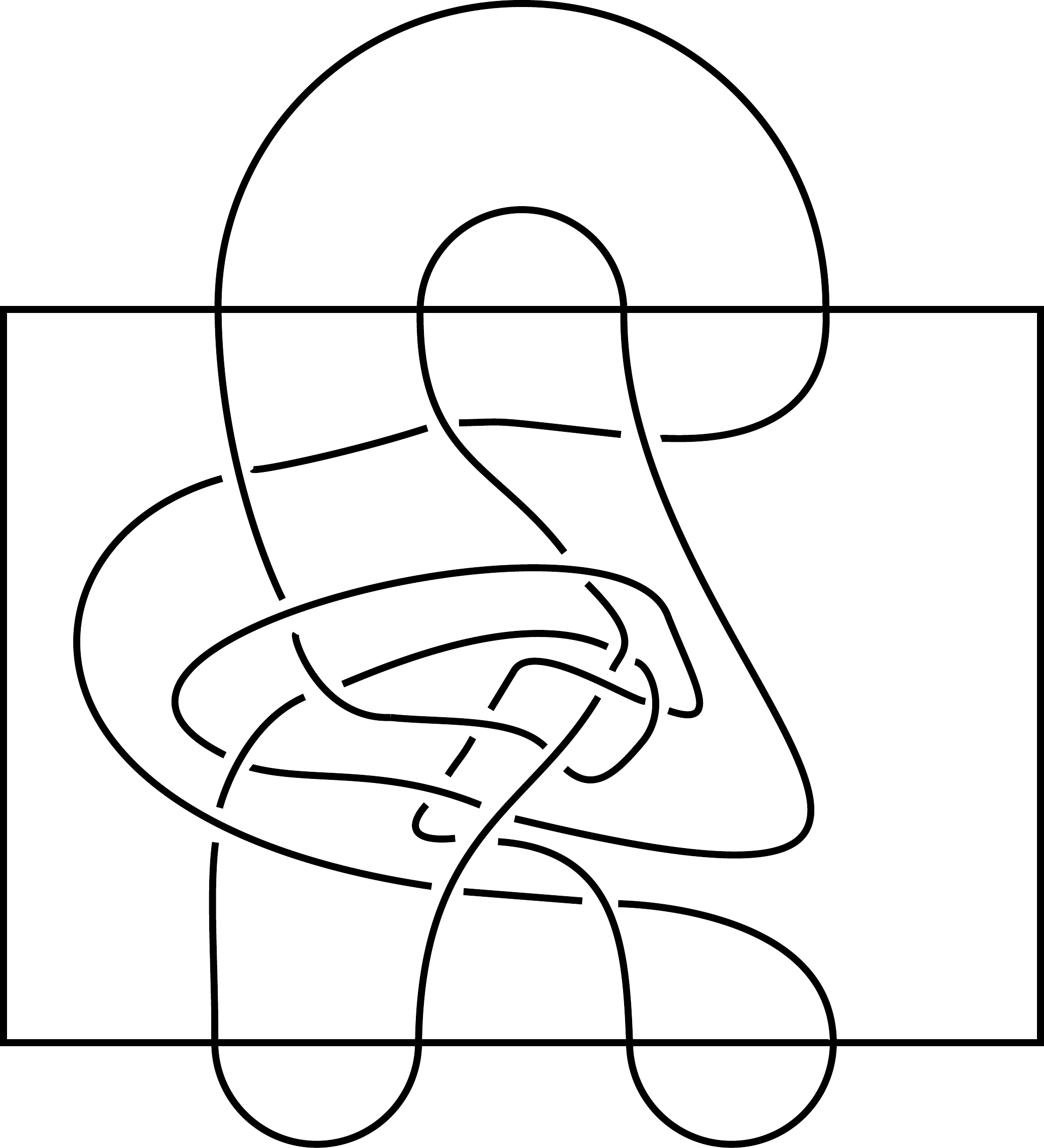}}%
  \end{picture}%
\endgroup%

%% file: diagrams/10_129.pdf_tex
\begingroup%
  \makeatletter%
  \providecommand\color[2][]{%
    \errmessage{(Inkscape) Color is used for the text in Inkscape, but the package 'color.sty' is not loaded}%
    \renewcommand\color[2][]{}%
  }%
  \providecommand\transparent[1]{%
    \errmessage{(Inkscape) Transparency is used (non-zero) for the text in Inkscape, but the package 'transparent.sty' is not loaded}%
    \renewcommand\transparent[1]{}%
  }%
  \providecommand\rotatebox[2]{#2}%
  \ifx\svgwidth\undefined%
    \setlength{\unitlength}{599.275592bp}%
    \ifx\svgscale\undefined%
      \relax%
    \else%
      \setlength{\unitlength}{\unitlength * \real{\svgscale}}%
    \fi%
  \else%
    \setlength{\unitlength}{\svgwidth}%
  \fi%
  \global\let\svgwidth\undefined%
  \global\let\svgscale\undefined%
  \makeatother%
  \begin{picture}(1,1.09985008)%
    \put(0,0){\includegraphics[width=\unitlength,page=1]{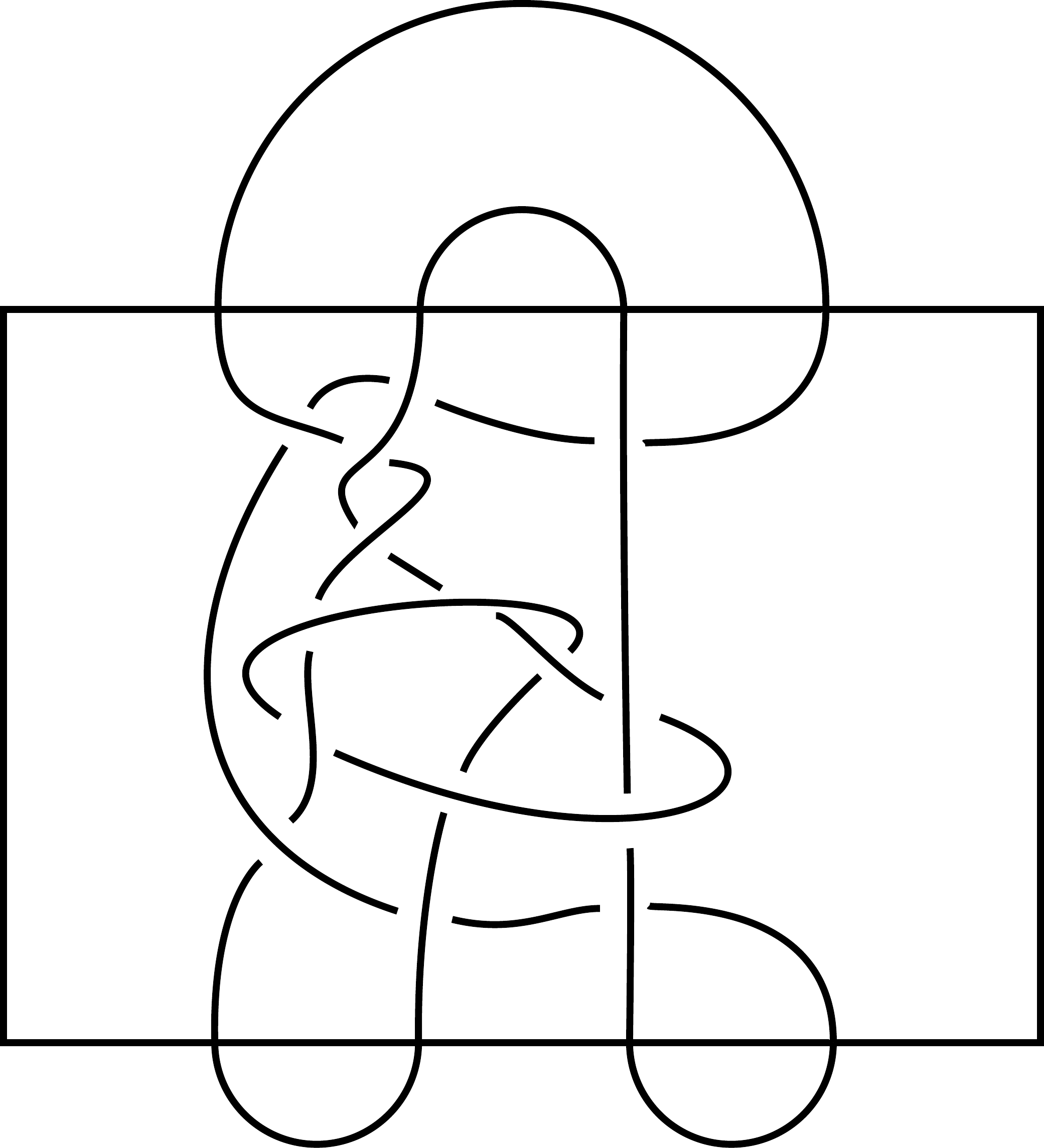}}%
  \end{picture}%
\endgroup%

%% file: diagrams/10_137.pdf_tex
\begingroup%
  \makeatletter%
  \providecommand\color[2][]{%
    \errmessage{(Inkscape) Color is used for the text in Inkscape, but the package 'color.sty' is not loaded}%
    \renewcommand\color[2][]{}%
  }%
  \providecommand\transparent[1]{%
    \errmessage{(Inkscape) Transparency is used (non-zero) for the text in Inkscape, but the package 'transparent.sty' is not loaded}%
    \renewcommand\transparent[1]{}%
  }%
  \providecommand\rotatebox[2]{#2}%
  \ifx\svgwidth\undefined%
    \setlength{\unitlength}{599.275592bp}%
    \ifx\svgscale\undefined%
      \relax%
    \else%
      \setlength{\unitlength}{\unitlength * \real{\svgscale}}%
    \fi%
  \else%
    \setlength{\unitlength}{\svgwidth}%
  \fi%
  \global\let\svgwidth\undefined%
  \global\let\svgscale\undefined%
  \makeatother%
  \begin{picture}(1,1.09985008)%
    \put(0,0){\includegraphics[width=\unitlength,page=1]{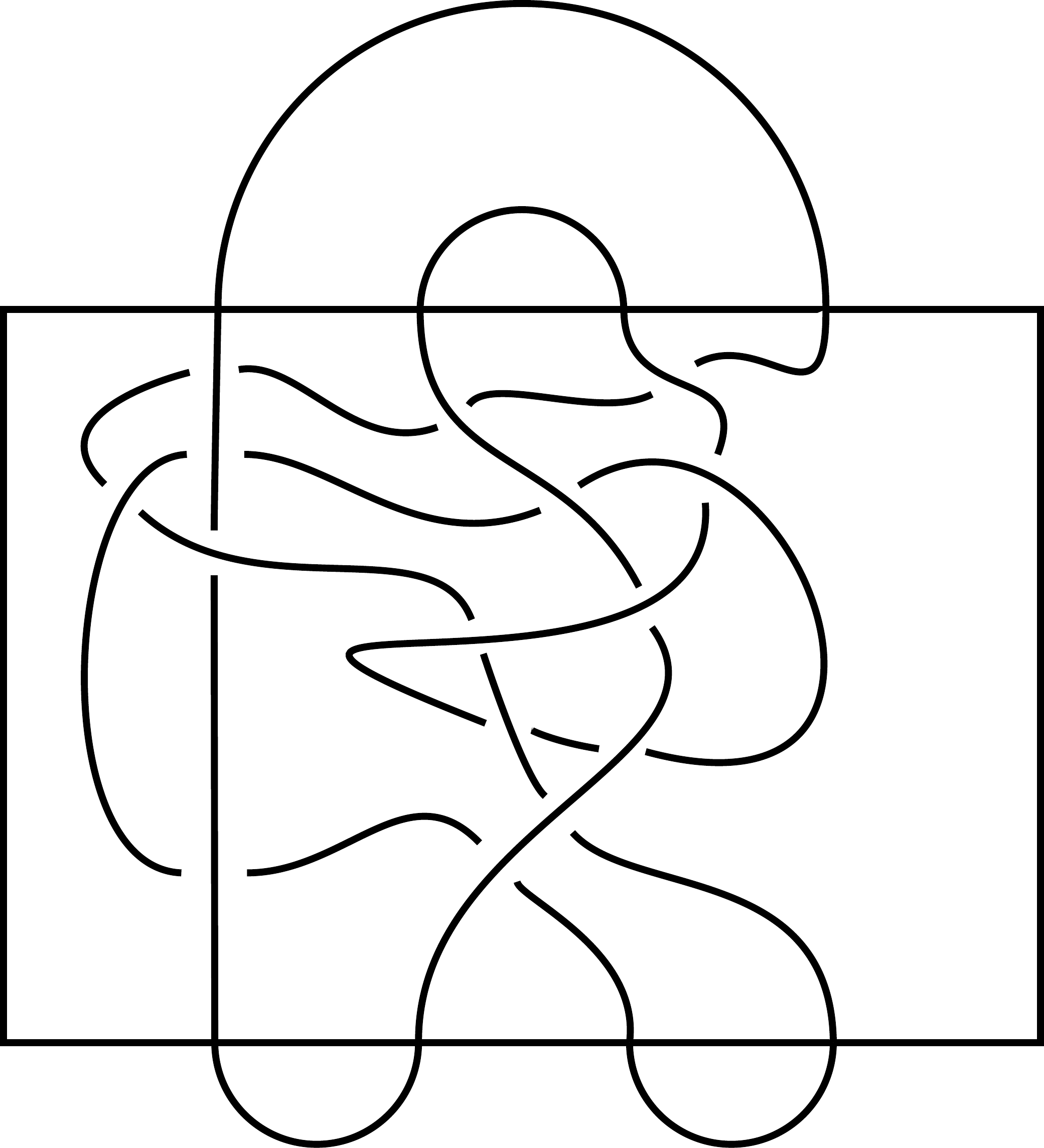}}%
  \end{picture}%
\endgroup%

%% file: diagrams/10_140.pdf_tex
\begingroup%
  \makeatletter%
  \providecommand\color[2][]{%
    \errmessage{(Inkscape) Color is used for the text in Inkscape, but the package 'color.sty' is not loaded}%
    \renewcommand\color[2][]{}%
  }%
  \providecommand\transparent[1]{%
    \errmessage{(Inkscape) Transparency is used (non-zero) for the text in Inkscape, but the package 'transparent.sty' is not loaded}%
    \renewcommand\transparent[1]{}%
  }%
  \providecommand\rotatebox[2]{#2}%
  \ifx\svgwidth\undefined%
    \setlength{\unitlength}{599.275592bp}%
    \ifx\svgscale\undefined%
      \relax%
    \else%
      \setlength{\unitlength}{\unitlength * \real{\svgscale}}%
    \fi%
  \else%
    \setlength{\unitlength}{\svgwidth}%
  \fi%
  \global\let\svgwidth\undefined%
  \global\let\svgscale\undefined%
  \makeatother%
  \begin{picture}(1,1.09985008)%
    \put(0,0){\includegraphics[width=\unitlength,page=1]{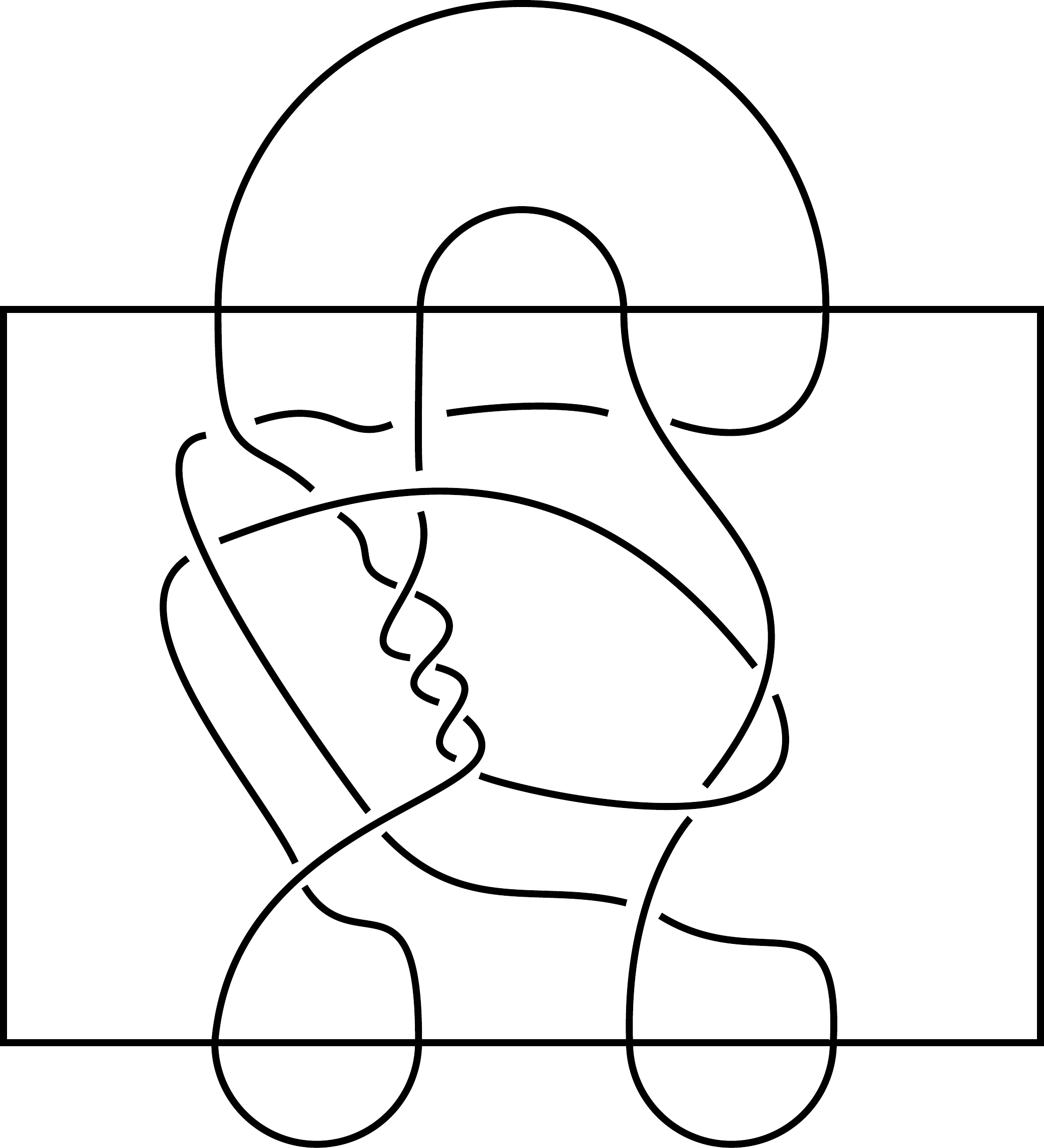}}%
  \end{picture}%
\endgroup%

%% file: diagrams/10_153.pdf_tex
\begingroup%
  \makeatletter%
  \providecommand\color[2][]{%
    \errmessage{(Inkscape) Color is used for the text in Inkscape, but the package 'color.sty' is not loaded}%
    \renewcommand\color[2][]{}%
  }%
  \providecommand\transparent[1]{%
    \errmessage{(Inkscape) Transparency is used (non-zero) for the text in Inkscape, but the package 'transparent.sty' is not loaded}%
    \renewcommand\transparent[1]{}%
  }%
  \providecommand\rotatebox[2]{#2}%
  \ifx\svgwidth\undefined%
    \setlength{\unitlength}{599.275592bp}%
    \ifx\svgscale\undefined%
      \relax%
    \else%
      \setlength{\unitlength}{\unitlength * \real{\svgscale}}%
    \fi%
  \else%
    \setlength{\unitlength}{\svgwidth}%
  \fi%
  \global\let\svgwidth\undefined%
  \global\let\svgscale\undefined%
  \makeatother%
  \begin{picture}(1,1.09985008)%
    \put(0,0){\includegraphics[width=\unitlength,page=1]{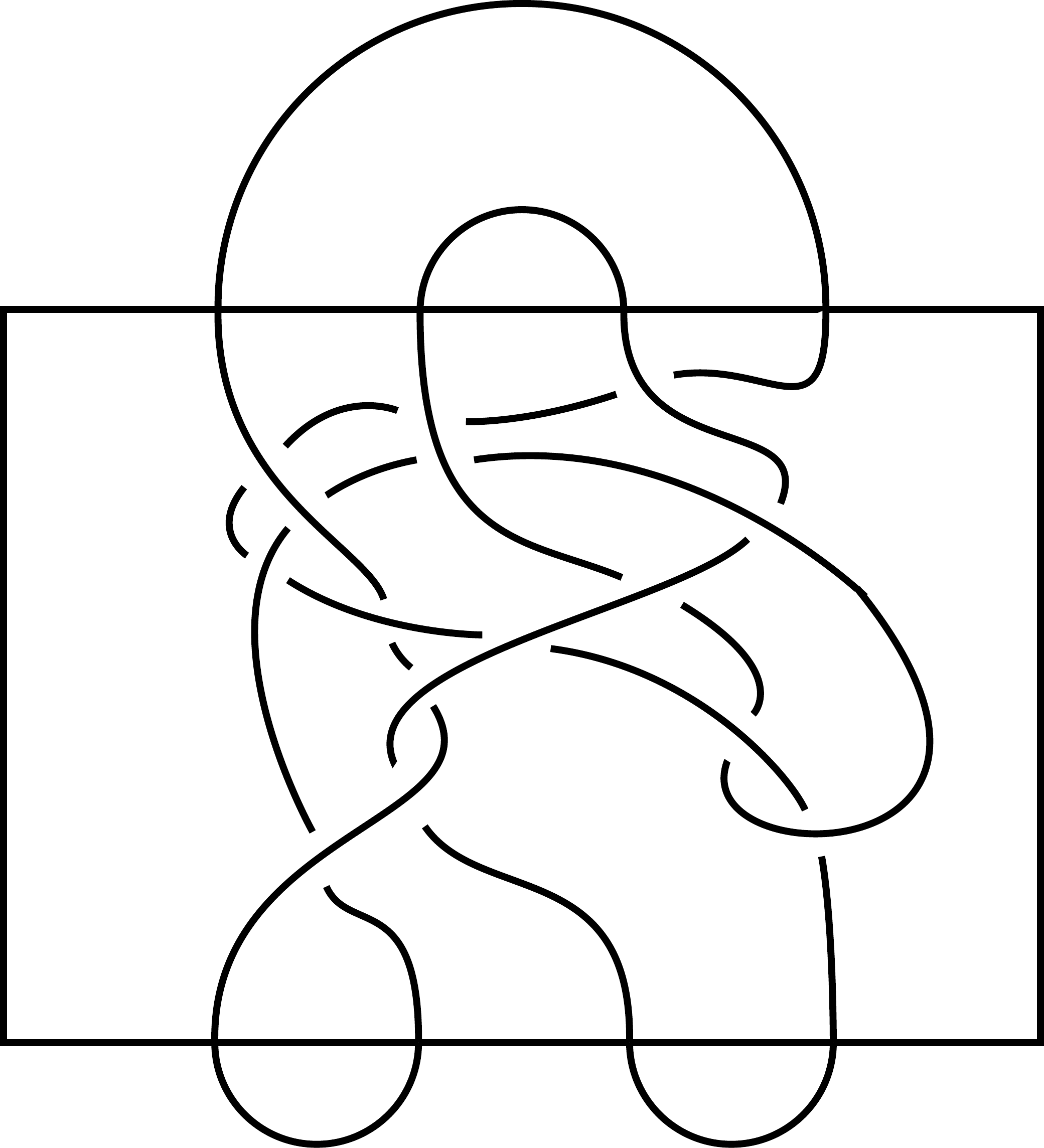}}%
  \end{picture}%
\endgroup%

%% file: diagrams/10_155.pdf_tex
\begingroup%
  \makeatletter%
  \providecommand\color[2][]{%
    \errmessage{(Inkscape) Color is used for the text in Inkscape, but the package 'color.sty' is not loaded}%
    \renewcommand\color[2][]{}%
  }%
  \providecommand\transparent[1]{%
    \errmessage{(Inkscape) Transparency is used (non-zero) for the text in Inkscape, but the package 'transparent.sty' is not loaded}%
    \renewcommand\transparent[1]{}%
  }%
  \providecommand\rotatebox[2]{#2}%
  \ifx\svgwidth\undefined%
    \setlength{\unitlength}{599.275592bp}%
    \ifx\svgscale\undefined%
      \relax%
    \else%
      \setlength{\unitlength}{\unitlength * \real{\svgscale}}%
    \fi%
  \else%
    \setlength{\unitlength}{\svgwidth}%
  \fi%
  \global\let\svgwidth\undefined%
  \global\let\svgscale\undefined%
  \makeatother%
  \begin{picture}(1,1.09985008)%
    \put(0,0){\includegraphics[width=\unitlength,page=1]{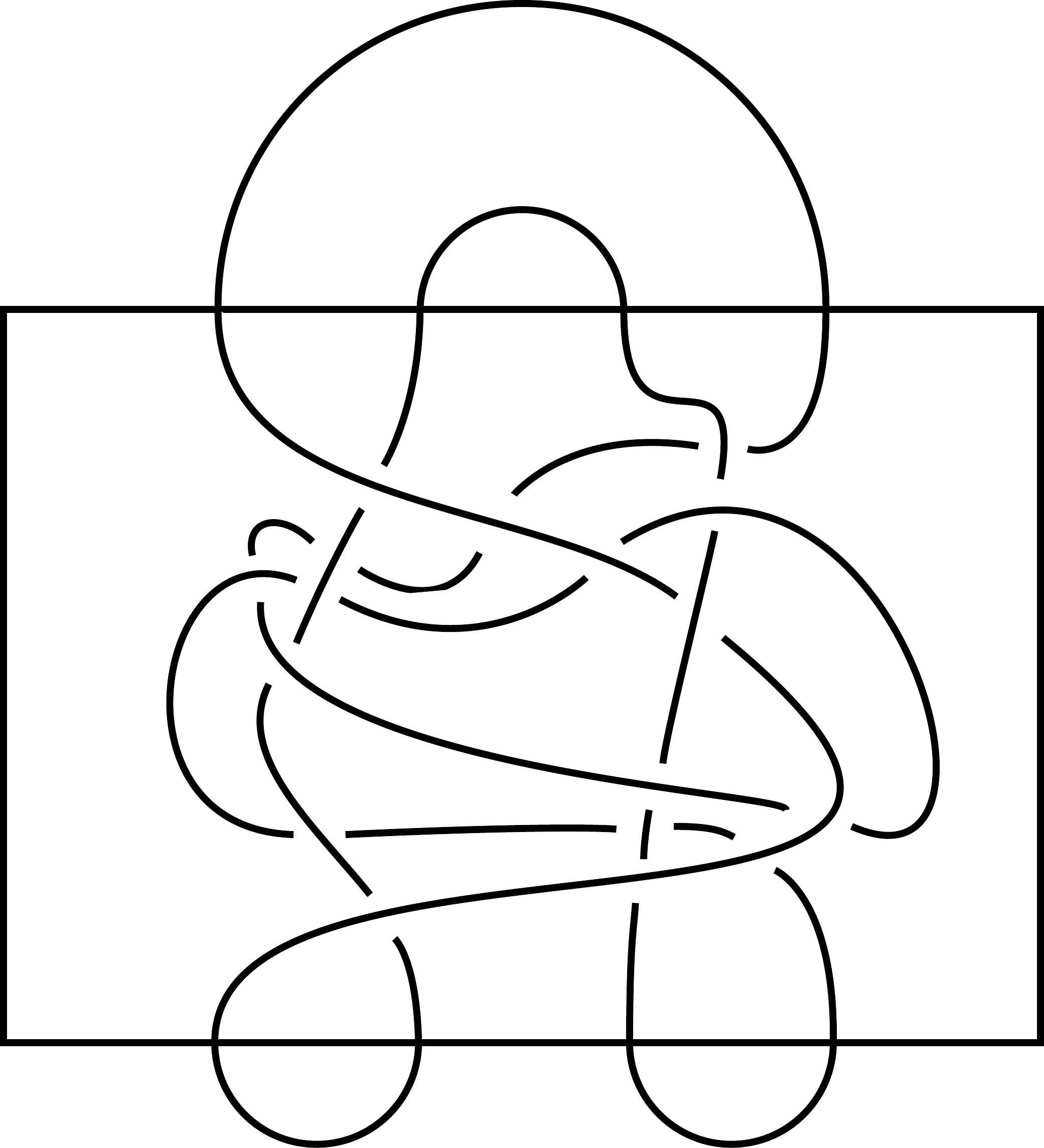}}%
  \end{picture}%
\endgroup%

%% file: ribbonknotdb.bbl
\begin{thebibliography}{}
\bibitem{knots}
Kawauchi, Akio.
\textit{Survey of Knot Theory.}
N.p.: Birkhauser, 2012. Print.
\bibitem{one}
Lamm, Christoph.
\textit{Symmetric union presentations for 2-bridge ribbon knots.}
\texttt{arXiv:math/0602395}
\bibitem{many}
Lamm, Christoph.
\textit{Symmetric unions and ribbon knots.}
Osaka J. Math. 37 (2000), no. 3, 537--550
\end{thebibliography}
